\documentclass{amsart}
\usepackage{hyperref}
\usepackage{graphicx}
\usepackage{amsmath,amssymb,amsthm}
\usepackage{color,subfigure}
\usepackage{algorithm,algorithmic}
\def\argmax{\operatornamewithlimits{arg\,max}}
\def\<#1,#2>{\langle#1,#2\rangle}
\newcommand{\mb}[1]{\mathbf{#1}}

\newcommand{\mrm}[1]{\text{\rm #1}}
\newcommand{\dd}{\,d}
\newcommand{\col}{\operatorname{col}}

\newcommand{\new}[1]{{\em #1}}

\newcommand{\R}{\mathbb{R}}
\newcommand{\rbar}{\overline{\R}}
\newcommand{\Z}{\mathbb{Z}}
\newcommand{\sK}{\mathcal{K}}

\newcommand{\im}{\mathrm{im}\,}
\newcommand{\projim}[1]{\Pi_{#1}}

\newcommand{\projimker}[2]{\Pi_{#1}^{#2}}
\newcommand{\rmax}{\mathbb{R}_{\max}}
\newcommand{\rmaxb}{\overline{\R}_{\max}}
\newcommand{\rminb}{\overline{\R}_{\min}}

\newcommand{\supp}{\mathop{\text{\Large$\vee$}}}
\newcommand{\inff}{\mathop{\text{\Large$\wedge$}}}
\newcommand{\maxx}{\max}

\newcommand{\conv}{\mathrm{Conv}}
\newcommand{\ri}{\mathrm{ri}}    
\newcommand{\dom}{\mathrm{dom}}

\newcommand{\comp}{\circ}

\newcommand{\bydef}{\stackrel{\mathrm{def}}{=}}
\newcommand{\set}[2]{\{#1\mid\,#2\}}
\newcommand{\lres}{\backslash}

\newcommand{\sh}{^{\sharp}}

\newcommand{\calu}{\mathcal{U}}
\newcommand{\calv}{\mathcal{V}}
\newcommand{\calw}{\mathcal{W}}
\newcommand{\calx}{\mathcal{X}}
\newcommand{\caly}{\mathcal{Y}}
\newcommand{\calz}{\mathcal{Z}}
\def\<#1,#2>{\langle#1\mid #2\rangle}
\newtheorem{thm}{Theorem}
\newtheorem{prop}[thm]{Proposition}
\newtheorem{lem}[thm]{Lemma}
\newtheorem{cor}[thm]{Corollary}
\theoremstyle{definition}
\newtheorem{exmp}[thm]{Example}
\newtheorem{rem}[thm]{Remark}
\newtheorem{defin}[thm]{Definition}

\title[The max-plus finite element method]{The max-plus finite element method for solving deterministic optimal control problems: basic properties and convergence analysis}
\author{Marianne Akian}
\author{St\'ephane Gaubert}
\address{INRIA, Domaine de Voluceau, 78153 Le Chesnay C\'edex, France}
\email{\{Marianne.Akian,Stephane.Gaubert,Asma.Lakhoua\}@inria.fr}
\author{Asma Lakhoua}
\address{Asma Lakhoua is also at ENIT-LAMSIN, BP 37, 1002 Tunis Le Belv\'ed\`ere, Tunisie}
\email{Asma.Lakhoua@lamsin.rnu.tn}
\date{March 27th, 2006}

\thanks{This work was supported by a fellowship of the IFC (Institut Fran\c{c}ais de Coop\'eration), by a fellowship of the AUF (Agence Universitaire de la Francophonie) and by the grant STIC-INRIA-Universit\'es tunisiennes \#I-04}
\begin{document}
\begin{abstract}
We introduce a max-plus analogue of the Petrov-Galerkin finite element method to solve finite horizon deterministic optimal control problems. The method relies on a max-plus variational formulation. We show that the error in the sup norm can be bounded from the difference between the value function and its projections on max-plus and min-plus semimodules, when the max-plus analogue of the stiffness matrix is exactly known. In general, the stiffness matrix must be approximated: this requires approximating the operation of the Lax-Oleinik semigroup on finite elements. We consider two approximations relying on the Hamiltonian. We derive a convergence result, in arbitrary dimension, showing that for a class of problems, the error estimate is of order $\delta+\Delta x(\delta)^{-1}$ or $\sqrt{\delta}+\Delta x(\delta)^{-1}$, depending on the choice of the approximation, where $\delta$ and $\Delta x$ are respectively the time and space discretization steps. We compare our method with another max-plus based discretization method previously introduced by Fleming and McEneaney. We give numerical examples in dimension 1 and 2. 
\end{abstract}
\maketitle
\section{Introduction}
We consider the optimal control problem:
\begin{subequations}
\label{problemP}
\begin{align}
\label{p1}
\mrm{maximize }
\int_0^T \ell(\mathbf{x}(s),\mb u(s))\dd s+\phi(\mb x(T))
\end{align}
over the set of trajectories $(\mb x(\cdot),\mb u(\cdot))$ satisfying
\begin{align}
\label{p2}
\dot{\mb x}(s)&=f(\mb x(s),\mb u(s)),\quad \mb x(s)\in X,\quad \mb u(s)\in U \enspace,
\end{align}
for all $0\leq s\leq T$ and
\begin{align}
\label{p3}
\mb x(0)=x \enspace.
\end{align}
\end{subequations}
Here, the \new{state space} $X$ is a subset of $\R^n$,
the set of \new{control values} $U$ is a subset of $\R^m$,
the \new{horizon} $T>0$ and the \new{initial condition} $x\in X$ are given,
we assume that the map $\mb u(\cdot)$ is measurable, 
and that the map $\mb x(\cdot)$ is absolutely continuous.
We also assume that the \new{instantaneous reward}
or \new{Lagrangian} $\ell:X\times U \to \R$,
and the \new{dynamics} $f:X\times U \to \R^n$, are
sufficiently regular maps, and that the \new{terminal reward}
$\phi$ is a map $X\to \R\cup\{-\infty\}$.

We are interested in the numerical computation of the \new{value function} $v$ which associates to any $(x,t)\in X\times [0,T]$ the supremum $v(x,t)$ of $\int_0^t \ell(\mb x(s),\mb u(s))\dd s+\phi(\mb x(t))$,
under the constraints~\eqref{p2}, for $0\leq s \leq t$ and~\eqref{p3}. It is known that, under certain regularity assumptions, $v$ is solution of the
Hamilton-Jacobi equation
\begin{subequations}\label{HJ}
\begin{gather}
-\frac{\partial v}{\partial t}+H(x,\frac{\partial v}{\partial x})=0, \quad (x,t)
\in X \times (0,T] \enspace,\label{HJ1}\end{gather}
with initial condition:
\begin{gather}
v(x,0)=\phi(x), \quad  x \in X \enspace,
\end{gather}
\end{subequations}
where $H(x,p)=\sup_{u\in U}\ell(x,u)+p\cdot f(x,u)$ is the \new{Hamiltonian}
of the problem (see for instance~\cite{lions,soner,barles}).

Several techniques have been proposed in the litterature to solve this problem. We mention for example finite difference schemes and the method of the vanishing viscosity~\cite{crandall-lions}, the anti-diffusive schemes for advection~\cite{zidani-bokanowski}, the finite elements approach~\cite{gonzalez-rofman} (in the case of the stopping time problem), the so-called discrete dynamic programming method or semi-lagrangian method~\cite{capuzzodolcetta}, \cite{capuzzodolcetta-ishii}, \cite{falcone}, \cite{falcone-ferretti}, \cite{falcone-giorgi}, \cite{carlini-falcone-ferretti}, the Markov chain approximations~\cite{boue-dupuis}. Other schemes have been obtained by integration from the essentially nonoscillatory (ENO) schemes for the hyperbolic conservation laws (see for instance~\cite{osher-shu}). Recently, max-plus methods have been proposed to solve first-order Hamilton-Jacobi equations~\cite{mceneaney-hortona}, \cite{mceneaney-hortonb}, \cite{mceneaney}, \cite{mceneaney02}, \cite{mceneaney03}, \cite{mceneaney-collins}, \cite{mceneaney04}.

Recall that the \new{max-plus semiring}, $\rmax$, is the set
$\R\cup\{-\infty\}$, equipped
with the addition $a\oplus b=\max(a,b)$ and the multiplication
$a\otimes b=a+b$. In the sequel, let $S^t$ denote the \new{evolution semigroup} of~\eqref{HJ}, or Lax-Oleinik semigroup, which associates
to any map $\phi$ the function $v^t:=v(\cdot,t)$, where $v$ is the value
function of the optimal control problem~\eqref{problemP}. Maslov~\cite{maslov73} observed that the semigroup $S^t$ is \new{max-plus linear}, meaning that for all maps $f,g$ from $X$ to $\rmax$, 
and for all $\lambda\in\rmax$, we have 
\begin{align*}
S^t(f\oplus g)&=S^tf\oplus S^tg \enspace ,\\
S^t(\lambda f)&=\lambda S^tf \enspace ,
\end{align*}
where $f\oplus g$ denotes the map $x\mapsto f(x)\oplus g(x)$,
and $\lambda f$ denotes the map $x\mapsto \lambda \otimes f(x)$.
Linear operators over max-plus type semirings have been widely studied,
see for instance~\cite{cuning,maslov92,baccelli,kolokoltsov,gondran-minoux}, see also~\cite{fathi}.

In~\cite{mceneaney}, Fleming and McEneaney introduced
a max-plus based discretization method to solve a subclass
of Hamilton-Jacobi equations
(with a Lagrangian $\ell$ quadratic with respect to $u$, and a dynamics
$f$ affine with respect to $u$). They use the max-plus linearity of the semigroup $S^t$ to approximate the value function $v^t$ by a function $v_h^t$ of the form:
\begin{equation}\label{v_form}
v_h^t=\sup_{1\leq i\leq p}\{\lambda_i^t+w_i\}\enspace ,
\end{equation}
where $\{w_i\}_{1\leq i\leq p}$ is a given family of functions (a max-plus ``basis'') and $\{\lambda_i^t\}_{1\leq i\leq p}$ is a family of scalars (the ``coefficients'' of $v_h^t$ on the max-plus ``basis''), which must be determined. They proposed a discretization scheme in which $\lambda^t$ is computed inductively by applying a max-plus linear operator to $\lambda^{t-\delta}$, where $\delta$ is the time discretization step. Thus, their scheme can be interpreted as the dynamic programming equation of a discrete control problem.

In this paper, we introduce a max-plus analogue of the finite element method, the ``MFEM'', to solve the deterministic optimal control problem~\eqref{problemP}. We still look for an approximation $v_h^t$ of the form~\eqref{v_form}. However, to determine the ``coefficients'' $\lambda_i^t$, we use a max-plus analogue of the notion of variational formulation,  which originates from the notion of generalized solution of Hamilton-Jacobi equations of Maslov and Kolokoltsov~\cite{kolokltsovmaslov88}, \cite[Section 3.2]{kolokoltsov}. We choose a family $\{z_j\}_{1\leq j\leq q}$ of test functions and define inductively $v_h^t$ to be the maximal function of the form~\eqref{v_form} satisfying
\begin{equation}
\<v_h^t,z_j>\leq\<S^\delta v_h^{t-\delta},z_j> \qquad \forall 1\leq j\leq q\enspace ,
\end{equation}
where $\<\cdot,\cdot>$ denotes the max-plus scalar product (see Section~\ref{method} for details). We show that the corresponding vector of coefficients $\lambda^t$ can be obtained by applying to $\lambda^{t-\delta}$ a nonlinear operator, which can be interpreted as the dynamic programming operator of a deterministic zero-sum two players game, with finite action and state spaces. The state space of the game corresponds to the set of finite elements. To each test function corresponds one possible action of the first player, and to each finite element corresponds one possible action of the second player, see Remark~\ref{jeu}. 

One interest of the MFEM is to provide, as in the case of the classical finite element method, a systematic way to compute error estimates, which can be interpreted geometrically as ``projection'' errors. In the classical finite element method, orthogonal projectors with respect to the energy norm must be used. In the max-plus case, projectors on semimodules must be used (note that these projectors minimize an additive analogue of Hilbert projective metric~\cite{ilade}). 

We shall see that when the value function is nonsmooth, the space of test functions must be different from the space in which the solution is represented, so that our discretization is indeed a max-plus analogue of the Petrov-Galerkin finite element method. A convenient choice of finite elements and test functions include quadratic functions (also considered by Fleming and McEneaney~\cite{mceneaney}) and norm-like functions, see Section~\ref{section_error}.

In the MFEM, we need to compute the value of the max-plus scalar product $\<z,S^\delta w>$ for each finite element $w$ and each test function $z$. In some special cases, $\<z,S^\delta w>$ can be computed analytically. In general, we need to approximate this scalar product. Here we consider the approximation $S^\delta w(x)=w(x)+\delta H(x,\nabla w(x))$, for $x\in X$, which is also used in~\cite{mceneaney-hortonb}. Our main result, Theorem~\ref{th-main}, provides for the resulting discretization of the value function an error estimate of order $\delta+\Delta x(\delta)^{-1}$, where $\Delta x$ is the ``space discretization step'', under classical assumptions on the control problem and the additionnal assumption that the value function $v^t$ is semiconvex for all $t\in[0,T]$. This is comparable with the order obtained in the simplest dicrete dynamic programming method, see~\cite{capuzzodolcetta-ishii}, \cite{falcone}, \cite{capuzzodolcetta-falcone}. To avoid solving a difficult (nonconvex) optimization problem, we propose a further approximation of the max-plus scalar product $\<z,S^\delta w>$, for which we obtain an error estimate of order $\sqrt\delta+\Delta x(\delta)^{-1}$, which is yet comparable to the order of the existing discretization methods~\cite{capuzzodolcetta-ishii}, \cite{falcone}, \cite{capuzzodolcetta-falcone}, \cite{crandall-lions}.

Note that the discretization grid need not be regular: in Theorem~\ref{th-main}, $\Delta x$ is defined for an arbitrary grid in term of Voronoi tesselations.

The paper is organised as follows. In Section~\ref{prelim}, we recall some basic tools and notions: residuation, semimodules and projection. In Section~\ref{method}, we present the formulation of the max-plus finite element method. 
In Section~\ref{compare_Mc} we compare our method with the method proposed by Fleming and McEneaney in~\cite{mceneaney}. In Section~\ref{section_error}, we state an error estimate and we give the main convergence theorem. Finally, in Section~\ref{Numerical results}, we illustrate the method by numerical examples
in dimension $1$ and $2$. Preliminary results of this paper appeared in~\cite{MTNS}.

\textbf{Acknowledgment:} We thank Henda El Fekih for advices and suggestions all along the development of the present work.

\section{Preliminaries on residuation and projections over semimodules}\label{prelim}
In this section we recall some classical residuation results (see
for example \cite{marie-louise}, \cite{birkhoff}, \cite{blyth},
\cite{baccelli}), and their application to linear maps on idempotent
semimodules (see~\cite{litvinov,ilade}). We also review some
results of \cite{wodes,ilade} concerning projectors over
semimodules. Other results on projectors over semimodules appeared in~\cite{gondran96a,gondran-minoux}.
\subsection{Residuation, semimodules, and linear maps}
If $(S,\leq)$ and $(T,\leq)$ are (partially) ordered sets, we say that a map
$f:S\to T$ is \new{monotone} if $s \leq s' \implies f(s)\leq
f(s')$. We say that $f$ is \new{residuated} if there exists a map
$f\sh: T\to S$ such that
\begin{equation}\label{e-def-res}
f(s) \leq t \iff s\leq f\sh(t) \enspace .
\end{equation}
The map $f$ is residuated if, and only if, for all $t\in T$,
$\set{s\in S}{f(s)\leq t}$ has a maximum element in $S$. Then,
 \begin{align*}
f\sh(t)&=\max\set{s\in S}{f(s)\leq t},\quad \forall t\in T \enspace.
\end{align*}
Moreover, in that case, we have 
\begin{align}\label{ffdf}
&f\comp f\sh\comp f=f\text{ and } f\sh\comp f\comp f\sh=f\sh\enspace.
\end{align}
In the sequel, we shall consider situations where $S$ (or $T$) is equipped with an idempotent monoid law $\oplus$ (\new{idempotent} means that $a\oplus a =a$). Then the \new{natural order} on $S$ is defined by $a\leq b \iff a\oplus b=b$. The supremum law for the natural order, which is denoted by $\supp$, coincides with $\oplus$ and the infimum law for the natural order, when it exists, will be denoted by $\inff$. We say that $S$ is \new{complete} as a naturally ordered set if any subset of $S$ has a least upper bound for the natural order.

If $\sK$ is an idempotent semiring, i.e.\, a semiring whose addition is idempotent, we say that the semiring $\sK$ is \new{complete} if it is complete as a naturally ordered set, and if the left and right multiplications: $\sK\to \sK$, $x\mapsto ax$ and $x\mapsto xa$, are residuated. Here and in the sequel, semiring multiplication is denoted by concatenation.

The max-plus semiring, $\rmax$,
is an idempotent semiring. It is not complete, but it can be embedded
in the complete idempotent semiring  $\rmaxb$ obtained 
by adjoining $+\infty$ to $\rmax$, with the convention that
$-\infty$ is absorbing for the multiplication.
The map $x\mapsto -x$ from $\rbar$ to itself yields an isomorphism
from $\rmaxb$ to the complete idempotent semiring $\rminb$,
obtained by replacing $\max$ by $\min$ and by exchanging the roles
of $+\infty$ and $-\infty$ in the definition of $\rmaxb$.

Semimodules over semirings are defined like modules over rings,
mutatis mutandis, see~\cite{litvinov,ilade}.
When $\sK$ is a complete idempotent semiring, we say that a (right)
$\sK$-semimodule $\calx$ is \new{complete} if it is complete as an idempotent monoid, and if, for all $u\in \calx$ and $\lambda\in \sK$,
the right and left multiplications, $R^\calx_{\lambda}:\;\calx\to \calx$,
$v\mapsto v\lambda$ and 
$L^\calx_{u}:\;\sK\to \calx$, $\mu\mapsto u\mu$, are residuated (for the natural order).
In a complete semimodule $\calx$, we define, for all $u,v\in \calx$,
\begin{align*}
  u\lres v &\bydef (L_u^\calx)\sh(v) = \max\set{\lambda\in \sK}{u\lambda \leq v} \enspace .
\end{align*}

We shall use \new{semimodules of functions}: when $X$ is a set and $\sK$ is a complete idempotent semiring, the set of functions $\sK^X$ is a complete $\sK$-semimodule for the componentwise addition $(u,v)\mapsto u\oplus v$ (defined by $(u\oplus v)(x)= u(x)\oplus v(x)$), and the componentwise multiplication $(\lambda,u)\mapsto u \lambda$ (defined by $(u\lambda)(x)= u(x)\lambda$).

If $\sK$ is an idempotent semiring, and if $\calx$ and $\caly$ are $\sK$-semimodules, we say that a map $A:\calx\to \caly$ is \new{linear}, or is a \new{linear operator}, if for all $u, v\in\calx$ and $\lambda, \mu\in\sK$, $A(u\lambda\oplus v\mu)=A(u)\lambda\oplus A(v)\mu$. 
Then, as in classical algebra, we use the
notation $Au$ instead of $A(u)$. 
When $A$ is residuated and
$v\in \caly$, we use the notation $A\backslash v$ or $A\sh v$ instead of
$A\sh (v)$. We denote by $L(\calx,\caly)$ the set of linear operators from $\calx$ to $\caly$. If $\sK$ is a complete idempotent semiring, if $\calx$, $\caly$, $\calz$ are complete $\sK$-semimodules, and if $A\in L(\calx,\caly)$ is residuated, then $L(\calx,\caly)$ and $L(\calx,\calz)$ are complete $\sK$-semimodules and the map $L_A:L(\calx,\caly)\to L(\calx,\calz)$, $B\mapsto A\comp B$, is residuated and we set $A\backslash C:=(L_A)\sh(C)$, for all $C\in L(\calx,\calz)$.

If $X$ and $Y$ are two sets, $\sK$
is a complete idempotent semiring, and $a\in \sK^{X\times Y}$,
we construct the linear operator $A$ from $\sK^Y$ to $\sK^X$ which associates
to any $u\in \sK^Y$ the function $Au\in \sK^X$ such that
$Au(x)=\supp_{y\in Y} a(x,y)u(y)$.
We say that $A$ is the \new{kernel operator} with \new{kernel}
or \new{matrix} $a$.
We shall often use the same notation $A$ for the operator and the kernel.
As is well known (see for instance~\cite{baccelli}), the kernel
operator $A$ is residuated, and
\[
(A\backslash v)(y)=\inff_{x\in X}A(x,y)\backslash v(x)\enspace .
\]
In particular, when $\sK=\rmaxb$, we have
\begin{align}
\label{e-conv}
(A\backslash v)(y)=\inf_{x\in X}(-A(x,y)+ v(x))= [- A^* (-v)](y)
\enspace,
\end{align}
where $A^*$ denotes the \new{transposed operator} $\sK^X\to \sK^Y$, 
which is associated to the kernel $A^*(y,x)=A(x,y)$. 
(In~\eqref{e-conv}, we use the convention that $+\infty$ is
absorbing for addition.)

\subsection{Projectors on semimodules}
Let $\sK$ be a complete idempotent semiring and $\calv$ denote a \new{complete subsemimodule} of a complete semimodule
$\calx$, i.e.\ a subset of $\calx$
that is stable by arbitrary sups and by the action of scalars.
We call \new{canonical projector} on $\calv$ the map
\begin{equation}\label{projecteur}
P_\calv: \calx\to \calx, \quad u\mapsto P_\calv(u) = \maxx\set{v\in \calv}{v\leq u}\enspace . 
\end{equation}
Let $W$ denote a \new{generating family} of a complete
subsemimodule $\calv$, which means that
any element $v\in \calv$ can be written as
$v=\supp\set{w\lambda_w}{w\in W}$, for some $\lambda_w\in\sK$.
It is known that
\[
P_\calv(u) =  \supp_{w\in W} w (w\lres u) 
\]
(see for instance~\cite{ilade}).
If $B:\calu\to\calx$ is a residuated linear operator, then when $\calu$ and $\calx$ are complete semimodules over $\sK$, the image $\im B$
of $B$ is a complete subsemimodule of $\calx$, and
\begin{equation}\label{PimB}
P_{\im B}=B\comp B\sh\enspace .
\end{equation}
The max-plus finite element methods relies
on the notion of projection on an image, parallel
to a kernel, which was introduced by Cohen, 
the second author, and Quadrat, in~\cite{wodes}.
The following theorem, of which Proposition~\ref{vhdelta} below
is an immediate corollary, is a variation
on the results of~\cite[Section~6]{wodes}.
\begin{thm}[Projection on an image parallel to a kernel]\label{piBC}
Let $\calu$, $\calx$ and $\caly$ be complete semimodules over $\sK$. Let $B:\calu\to\calx$ and $C:\calx\to\caly$ be two residuated linear
operators over $\sK$. Let $\projimker{B}{C}=B\comp(C\comp B)\sh\comp C$. We have
$\projimker{B}{C}=\projimker{B}{}\comp \projimker{}{C}$, where
$\projimker{B}{}=B\comp B\sh$ and $\projimker{}{C}=C\sh\comp
C$. Moreover, $\projimker{B}{C}$ is a projector, meaning that $(\projimker{B}{C})^2=\projimker{B}{C}$, and for all
$x\in\calx$:
\[
\projimker{B}{C}(x)=\maxx\set{y\in\im B}{Cy\leq Cx}\enspace .
\]
\end{thm}
\begin{proof}
The first assertion follows from $(C\comp B)\sh=B\sh\comp C\sh$. For the second assertion, we have
\begin{eqnarray*}
(\projimker{B}{C})^2 &=& \big(B\comp (C\comp B)\sh\comp C\big)\comp\big(B\comp(C\comp B)\sh\comp C\big)\\
&=& B\comp (C\comp B)\sh\comp C \quad (\mathrm{using}\quad\eqref{ffdf})\\
&=& \projimker{B}{C}\enspace .
\end{eqnarray*} 
To prove the last assertion, we use that $\projim{B}=P_{\im B}$ and~\eqref{e-def-res}, we deduce:
\begin{eqnarray*}
\projimker{B}{C}(x) &=& P_{\im B}\comp C\sh\comp C(x)\\
&=& \max\set{y\in\im B}{y\leq C\sh\comp C(x)}\\
&=& \max\set{y\in\im B}{Cy\leq Cx}\enspace .
\end{eqnarray*}
\end{proof}
The results of~\cite{wodes} characterize
the existence and uniqueness, for all $x\in X$,
of $y\in \im B$ such that $Cy=Cx$. In that case,
$y=\projimker{B}{C}(x)$.

When $\sK=\rmaxb$, and $C:\rmaxb^X\to\rmaxb^Y$
is a kernel operator,
$\projimker{}C=C\sh\comp C$ has an interpretation similar
to~\eqref{PimB}:
\[ 
\projimker{}{C}(v)=C\sh\comp C(v)=-P_{\im C^*}(-v)
=P^{-\im C^*} (v)\enspace,
\]
where $-\im C^*$ is thought of as a $\rminb$-subsemimodule
of $\rminb^X$ and $P^\calv$ denotes the projector on a $\rminb$-semimodule $\calv$, so that, 
\[
P^{-\im C^*} (v)= \min\set{w\in -\im C^*}{w\geq v} \enspace ,
\]
where $\leq$ denotes here the usual order
on $\rbar^X$.
When $B:\rmaxb^U\to\rmaxb^X$ is also a kernel
operator, we have
\[
\projimker BC=P_{\im B}\comp P^{-\im C^*} \enspace .
\]
This factorization will be instrumental in the geometrical
interpretation of the finite element algorithm.
\begin{exmp}
We take $U=\{1,\cdots,p\}$, $X=\R$ and $Y=\{1,\cdots,q\}$. Consider the linear operators $B:\rmaxb^U\to\rmaxb^X$ and $C:\rmaxb^X\to\rmaxb^Y$ such that
\[
B\lambda(x)=\sup_{1\leq i\leq p}\{-\frac c 2(x-\hat x_i)^2+\lambda_i\},\qquad \text{for all }\lambda\in\rmaxb^U\enspace ,
\]
and
\[
(Cf)_i=\sup_{x\in\R}\{-a|x-\hat y_i|+f(x)\}, \qquad \text{for all } f\in\rmaxb^X\enspace .
\]
The image of $B$, $\im B$, is the semimodule generated in the max-plus sense by the functions $x\mapsto -\frac c 2(x-\hat x_i)^2$, for $i=1,\cdots,p$. We have
\[
C\sh\mu(x)=\inf_{1\leq i\leq q}\{a|x-\hat y_i|+\mu_i\}, \qquad \text{for all }\mu\in\rmaxb^Y\enspace ,
\]
and the image of $C\sh$, which coincides with $-\im C^*$, is the semimodule generated in the min-plus sense by the functions $x\mapsto a|x-\hat y_i|$, for $i=1,\cdots,q$.\\
In figure~\ref{proj-dessus-fig}, we represent a function $v$ and its projection $P^{-\im C^*}(v)$ (in bold). In figure~\ref{proj-dessous-fig}, we represent (in bold) the projection $P_{\im B}(P^{-\im C^*}(v))=\projimker{B}{C}(v)$.
\begin{figure}[htbp]
\centering
\mbox{\subfigure[]{\label{proj-dessus-fig}{\input{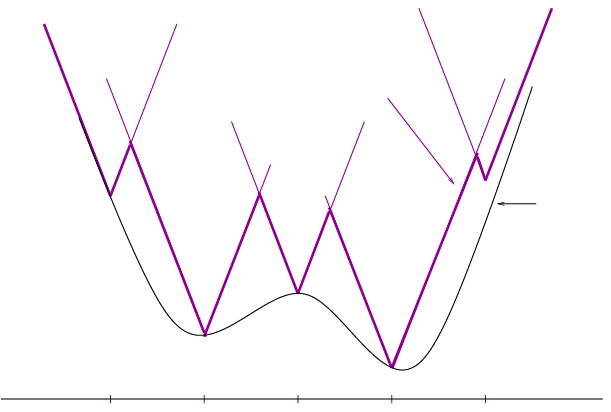}}}\quad
      \subfigure[]{\label{proj-dessous-fig}{\input{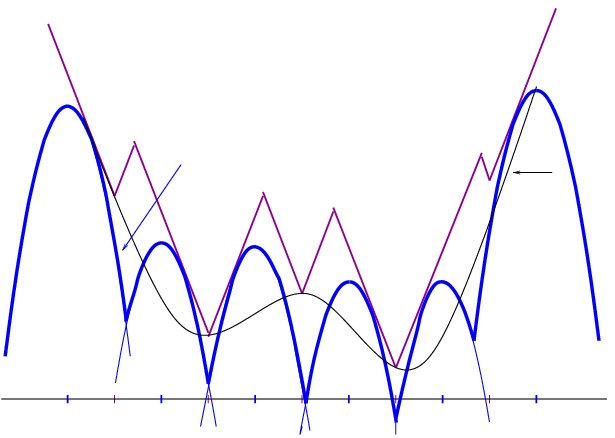}}}}
\caption{Example illustrating max-plus and min-plus projectors}
\label{example-projectors}
\end{figure}
\end{exmp}
\section{The max-plus finite element method}\label{method}
\subsection{Max-plus variational formulation}
We now describe the max-plus finite element method to solve Problem~(\ref{problemP}).
Let $\calv$ be a complete semimodule of functions from $X$ to $\rmaxb$. Let $S^t:\calv\to\calv$ and $v^t$ be defined as in the introduction.
Using the semigroup property $S^{t+t'}=S^t\circ S^{t'}$, for $t,t'>0$, 
we get:
\begin{equation}\label{exact}
\begin{array}{cc}
v^{t+\delta}=S^{\delta} v^t, & t=0,\delta, \cdots, T-\delta
\end{array}
\end{equation}
with $v^0=\phi$ and $\delta=\frac{T}{N}$, for some
positive integer $N$. Let $\calw\subset\calv$ be
a complete $\rmaxb$-semimodule of functions from $X$ to $\rmaxb$ such that
for all $t\geq0$, $v^t\in\calw$. We choose a ``dual'' semimodule $\calz$ of ``test functions'' from $X$
to $\rmaxb$. Recall that the max-plus \new{scalar product}
is defined by 
\[
\<u,v>=\sup_{x\in X} u(x)+v(x)\enspace ,
\]
for all functions $u,v:X\to \rmaxb$. We replace \eqref{exact} by:
\begin{equation}\label{produitscalaire}
\<z,v^{t+\delta}>=\<z,S^{\delta}v^{t}> ,
\quad  \forall z\in\calz \enspace ,
\end{equation}
for $t=0,\delta,\ldots,T-\delta$, with $v^{\delta},\ldots,v^T\in\calw$.
Equation~\eqref{produitscalaire}
can be seen as the analogue of a \new{variational}
or \new{weak formulation}. Kolokoltsov and Maslov used this formulation
in \cite{kolokltsovmaslov88} and \cite[Section 3.2]{kolokoltsov}
to define a notion of generalized solution of Hamilton-Jacobi equations.
\subsection{Ideal max-plus finite element method}
We consider a semimodule $\calw_h\subset\calw$ generated by the family $\{w_i\}_{1\leq i\leq p}$. We call \new{finite elements}
the functions $w_i$. We approximate $v^t$ by
$v_h^t\in\calw_h$, that is:
\[
v^t\simeq v_h^t=\supp_{1\leq i\leq p} w_i \lambda^t_i\enspace,
\]
where $\lambda_i^t\in\rmaxb$. We also consider a semimodule
$\calz_h\subset\calz$ with generating family $\{z_j\}_{1\leq j\leq
  q}$. The functions $z_1,\cdots,z_q$ will act as test functions. We replace \eqref{produitscalaire} by
\begin{equation}\label{approxgeneral}
\<z,v_h^{t+\delta}>=\<z,S^{\delta}v_h^{t}> ,
\quad  \forall z\in\calz_h \enspace ,
\end{equation}
for $t=0,\delta,\ldots,T-\delta$, with $v_h^{\delta},\ldots,v_h^T\in\calw_h$. The function $v_h^0$ is a given approximation of $\phi$. Since $\calz_h$ is generated by $z_1,\ldots,z_q$, \eqref{approxgeneral} is equivalent to
\begin{equation}\label{pdtscalapproch}
\begin{array}{cc}
\<z_j,v_h^{t+\delta}>=\<z_j,S^{\delta}v_h^{t}>, & \forall 1\leq
j\leq q\enspace,\\
\end{array}
\end{equation}
for $t=0,\delta, \cdots, T-\delta$, with $v_h^t\in\calw_h$, $t=0,\delta,\cdots,T$.

Since Equation \eqref{pdtscalapproch} need not have a
solution, we look for its maximal subsolution,
i.e.\ the maximal solution $v_h^{t+\delta}\in \calw_h$ of
\begin{subequations}\label{inegal}
\begin{align}
\<z_j,v_h^{t+\delta}> \quad \leq \quad \<z_j,S^{\delta}v_h^{t}> \quad 
\forall 1\leq j\leq q\enspace. \label{infouegal}
\end{align}
We also take for the approximate value function $v_h^0$ at time $0$
the maximal solution $v_h^0\in \calw_h$ of
\begin{align}
v_h^{0}\leq v^0\enspace.\label{inegvh0}
\end{align}
\end{subequations}
Let us denote by $W_h$ the max-plus linear operator from
$\rmaxb^p$ to $\calw$ with matrix $W_h=\col(w_{i})_{1\leq i\leq p}$,
and by $Z_h^*$ the max-plus linear operator from $\calw$ to $\rmaxb^q$ 
whose transposed matrix is $Z_h=\col(z_j)_{1\leq j\leq q}$.
This means that 
$W_h\lambda=\supp_{1\leq i\leq p} w_i \lambda_i$ for all
$\lambda=(\lambda_i)_{i=1,\ldots, p}\in \rmaxb^p$, and
$(Z_h^* v)_j=\< z_j,v>$ for all $v\in\calw$ and $j=1,\ldots,q$.
Applying Theorem~\ref{piBC} to $B=W_h$ and $C=Z_h^*$ and noting that
$\calw_h=\im W_h$, we get:
\begin{cor}\label{vhdelta}
The maximal solution $v_h^{t+\delta}\in \calw_h$
of~\eqref{infouegal} is given by 
$v_h^{t+\delta}=S_{h}^{\delta}(v_h^t)$, where 
\begin{equation*}
S_{h}^{\delta}=\projimker{W_h}{Z_h^*}\comp 
 S^{\delta}\enspace.
\end{equation*}
\end{cor}
Note that $\projimker{W_h}{Z_h^*}=P_{\calw_h}\comp P^{-\calz_h}$.
The following proposition provides a recursive equation verified by the vector of coordinates of $v_h^t$.
\begin{prop}Let $v_h^{t}\in \calw_h$ be the  maximal solution 
of~\eqref{inegal}, for $t=0,\delta,\ldots, T$. 
Then,  for every $t=0,\delta,\ldots, T$,
there exists a maximal $\lambda^t\in\rmax^p$ such that
$v_h^t=W_h\lambda^t$, $t=0,\delta,\cdots,T$, which can be determined recursively from 
\begin{subequations}
\begin{align}\label{lambdat+dt}
\lambda^{t+\delta}=(Z_h^*W_h)\backslash (Z_h^*S^{\delta}W_h\lambda^{t})
\enspace,
\end{align}
for $t=0,\cdots,T-\delta$ with the initial condition:
\begin{align}
\lambda^{0}=W_h\backslash \phi\enspace.
\end{align}
\end{subequations}
\end{prop}
\begin{proof}
Since $v_h^t\in \calw_h$, $v_h^t=W_h\lambda^t$ for some $\lambda^t\in\rmaxb^p$ 
and the maximal $\lambda^t$ satisfying this condition is
$\lambda^t=W_h\sh (v_h^t)$, for all $t=0,\delta,\ldots, T$.
Since  $v_h^0$ is the maximal solution of~\eqref{inegvh0},
then by~\eqref{projecteur} 
and~\eqref{PimB}, $v_h^0=P_{\calw_h}(\phi)=W_h \comp W_h\sh (\phi)$, hence
$\lambda^0= W_h\sh\comp W_h\comp W_h\sh (\phi)=W_h\sh(\phi)$. 
Let $t=\delta,\ldots, T$. Using Proposition~\ref{vhdelta},  
Theorem~\ref{piBC},~\eqref{ffdf} and the property that
$(f\comp g)\sh=g\sh\comp f\sh$ for all residuated maps $f$ and $g$,
we get 
\begin{eqnarray*}
\lambda^{t+\delta} & = & W_h\sh \comp \projimker{W_h}{Z_h^*}
\comp S^{\delta}(W_h\lambda^{t} )\\
& = & W_h\sh\comp W_h\comp  W_h\sh \comp (Z_h^*)\sh\comp Z_h^*
\comp S^{\delta}(W_h\lambda^{t})\\
&= & W_h\sh \comp (Z_h^*)\sh\comp Z_h^*\comp S^{\delta} (W_h\lambda^{t})\\
& = & (Z_h^* W_h)\sh (Z_h^* S^{\delta}W_h\lambda^{t})\enspace.
\end{eqnarray*}
which yields~\eqref{lambdat+dt}.
\end{proof}

For $1\leq i\leq p$ and $1\leq j\leq q$, we define:
\begin{align}
(M_h)_{ji}&=\<z_j,w_i> \label{matrixA}\\
(K_h)_{ji}&=\<z_j,S^{\delta}w_i> \label{matrixB}
\\
&=\<(S^*)^{\delta} z_j, w_i>
\enspace,\label{matrixB'}
\end{align}
where $S^*$ is the \new{transposed semigroup} of $S$, which
is the evolution semigroup associated to the optimal
control problem~(\ref{problemP}) in which the sign of the dynamics is changed. The matrices $M_h$ and $K_h$ represent respectively the max-plus linear operators $Z_h^*W_h$ and $Z_h^*S^\delta W_h$. Equation~\eqref{lambdat+dt} may be written explicitly, for $1\leq i\leq p$, as
\[
\lambda^{t+\delta}_i = \min_{1\leq j\leq q}
\Big( -(M_{h})_{ji}+ \max_{1\leq k\leq p}\big( (K_h)_{jk} + \lambda^{t}_k\big) \Big)\enspace .
\]
\begin{rem}\label{jeu}
This recursion may be interpreted as the dynamic programming equation of a deterministic zero-sum two players game, with finite action and state spaces. Here the state space of the game is the finite set $\{1,\cdots,p\}$ (to each finite element corresponds a state of the game). To each test function corresponds one possible action $j\in\{1,\cdots,q\}$ of the first player, and to each finite element corresponds one possible action $k\in\{1,\cdots,p\}$ of the second player. Given these actions at the state $i\in\{1,\cdots,p\}$, the cost of the first player, which is the reward of the second player, is $-(M_h)_{ji}+(K_h)_{jk}$.
\end{rem}
The ideal max-plus finite element method can be summarized as follows:
\begin{algorithm}
\caption{Ideal max-plus finite element method}
\begin{algorithmic}[1]
\STATE Choose the finite elements $(w_i)_{1\leq i\leq p}$ and $(z_j)_{1\leq j\leq q}$. Choose the time discretization step $\delta=\frac{T}{N}$,
\STATE Compute the matrix $M_h$ by~\eqref{matrixA}
and the matrix $K_{h}$ by~\eqref{matrixB} or by~\eqref{matrixB'},
\STATE Compute $\lambda^{0}=W_h\backslash\phi$ and
 $v_h^{0}=W_h\lambda^{0}$.
\STATE For $t=\delta, 2\delta,\ldots,T$, compute $\lambda^{t}=M_h\backslash (K_{h}\lambda^{t-\delta})$ and $v_h^{t}=W_h\lambda^{t}$.
\end{algorithmic}
\end{algorithm}
\begin{rem}
Since $v_h^t\in\calw_h$ $\forall t=0,\cdots,T$, the dynamics of $v_h^t$can be written as a function of the matrices $M_h$ and $K_h$:
\begin{equation}\label{v_MK}
v_h^{t+\delta}=W_h\comp M_h\sh\comp K_h\comp W_h\sh(v_h^t)\enspace .
\end{equation}
\end{rem}
\subsection{Effective max-plus finite element method}
In order to implement the max-plus finite element method, we must specify how to compute the entries of the matrices $M_h$ and $K_h$ in~\eqref{matrixA}
and~\eqref{matrixB} or~\eqref{matrixB'}.\\
Computing $M_h$ from~\eqref{matrixA} is an optimization
problem, whose objective function is concave for natural choices of finite elements and test functions (see Section~\ref{section_error} below). This problem may be solved by standard optimization algorithms. Evaluating every scalar product $\<z,S^\delta w>$ leads to a new optimal control problem since
\begin{equation*}
\<z,S^\delta w>=\max z(\mb x(0))+\int_0^\delta\ell(\mb x(s),\mb u(s))ds+w(\mb x(\delta))\enspace ,
\end{equation*}
where the maximum is taken over the set of trajectories $\big(\mb x(\cdot),\mb u(\cdot)\big)$ satisfying~\eqref{p2}. This problem is simpler to approximate than Problem~(\ref{problemP}), because the horizon $\delta$ is small, and the functions $z$ and $w$ have a regularizing effect. \\
We first discuss the approximation of $S^{\delta}w$ for every finite element $w$.
The Hamilton-Jacobi equation~\eqref{HJ1} suggests to
approximate $S^{\delta}w$ by the function $[S^{\delta}w]_H $ such that
\begin{equation}\label{stilde}
[S^{\delta}w]_H (x)=w(x)+\delta H(x,\nabla w(x)),
\quad \mrm{for all } x\in X\enspace  .
\end{equation}
Let $[S^{\delta}W_h]_H$ denote the max-plus linear operator from
$\rmaxb^p$ to $\calw$ with matrix $[S^{\delta}W_h]_H =
\col([S^{\delta}w_i]_H )_{1\leq i\leq p}$, which means that 
\[
[S^{\delta}W_h]_H \lambda=\supp_{1\leq i \leq p}[S^{\delta}w_i]_H\lambda_i
\]
 for all $\lambda=(\lambda_i)_{1\leq i\leq p}\in \rmax^p$.
The above approximation of  $S^{\delta}w$ 
yields an approximation of the matrix $K_h$ 
by the matrix $K_{H,h}:= Z_h^* [S^{\delta}W_h]_H$, whose
entries are given, for $1\leq i \leq p$ and $1\leq j\leq q$, by:
\begin{eqnarray}
(K_{H,h})_{ji}&=&
\sup_{x\in X}(z_j(x)+w_i(x)+\delta H(x,\nabla w_i(x)))\enspace.
\label{appr_K}
\end{eqnarray}
Thus, computing $K_{H,h}$ requires to solve an optimization problem,
which is nothing but a perturbation of the optimization problem
associated to the computation of $M$. 
We may exploit this observation by replacing
$K_{H,h}$ by the matrix $\tilde K_{H,h}$ with entries
\begin{eqnarray}
(\tilde K_{H,h})_{ji}&=&\<z_j,w_i>+\delta\sup_{x\in\argmax\{z_j+w_i\}}H(x,\nabla w_i(x))\enspace,
\label{e-convenient}
\end{eqnarray}
for $1\leq i \leq p$ and $1\leq j\leq q$. Here, 
$\argmax\{z_j+w_i\}$ denotes the set of $x$ such that 
$z_j(x)+w_i(x)=\<z_j,w_i>$. When this set has only one element,
\eqref{e-convenient} yields a convenient approximation
of $K_h$.

Of course, $w_i$ must be differentiable for the approximation~\eqref{stilde}
to make sense. When $w_i$ is non-differentiable, but $z_j$ is differentiable, the dual formula~\eqref{matrixB'} suggests to approximate $(K_h)_{ji}$ by
\begin{eqnarray*}
\sup_{x\in X}(z_j(x)+\delta H(x,-\nabla z_j(x)) +w_i(x))\enspace.
\end{eqnarray*}
We may also use the dual 
formula of~\eqref{e-convenient}, where $\nabla w_i(x)$ is replaced by $-\nabla z_j(x)$.
\section{Comparison with the method of Fleming and McEneaney}\label{compare_Mc}
Fleming and McEneaney proposed a max-plus based
method~\cite{mceneaney}, which also uses a space $\calw_h$ generated
by finite elements, $w_1,\ldots,w_p$, together with the linear
formulation~\eqref{exact}.
Their method approaches the value function at time $t$, 
$v^t$, by $W_h \mu^t$, where $W_h=\col(w_i)_{1\leq i\leq p}$ as
above, and $\mu^t$ is defined inductively by
\begin{subequations}\label{algo-fm}
\begin{align}
\mu^0 &= W_h\lres \phi \\
\mu^{t+\delta} & =  \big(W_h\backslash (S^{\delta}W_h)\big)\mu^{t} \enspace,\label{mut+dt}
\end{align}\end{subequations}
for $t=0,\delta,\ldots, T-\delta$. This can be compared with the limit
case of our finite element method, in which
the space of test functions $\calz_h$ is the set of all functions. This limit case
corresponds to replacing $Z_h^*$ by the identity operator
in~\eqref{lambdat+dt}, so that 
\begin{align}
\label{e-limit}
\lambda^{t+\delta}=W_h\backslash (S^{\delta}W_h\lambda^t) \enspace .
\end{align}
\begin{prop}
Let $(\mu^t)$ be the sequence of vectors
defined by the algorithm of Fleming and McEneaney, \eqref{algo-fm};
let $(\lambda^t)$ be the sequence of vectors defined
by the max-plus finite element method, in the limit case~\eqref{e-limit};
and let $v^t$ denote the value function at time $t$. 
Then, 
\begin{equation}\label{comparaison}
W_h \mu^t \leq W_h \lambda^t \leq v^t \enspace ,
\quad \mrm{for } t=0,\delta, \ldots, T \enspace .
\end{equation}
\end{prop}
\begin{proof}
We first prove that $W_h\lambda^t\leq v^t$ for $t=0,\delta, \ldots, T$. This can be proved by induction. For $t=0$ we have $W_h\lambda^0\leq v^0$ by~\eqref{inegvh0}. We assume that $W_h\lambda^t\leq v^t$. Using~\eqref{e-limit}, we have
\begin{eqnarray*}
W_h\lambda^{t+\delta} & =& W_hW_h\sh S^\delta(W_h\lambda^t)\\
& =& \projimker{W_h}{}\big(S^\delta(W_h\lambda^t)\big)\enspace .
\end{eqnarray*}
Using the monotonicity of the semigroup $S^\delta$, we obtain
\begin{eqnarray*}
W_h\lambda^{t+\delta} & \leq & \projimker{W_h}{}\big(S^\delta v^t\big)\\
& \leq & S^\delta v^t\\
& = & v^{t+\delta}\enspace .
\end{eqnarray*}
The second inequality is also proved by induction. For $t=0$, we have $\mu^0=\lambda^0=W_h\backslash\Phi$. Suppose that $\mu^t\leq\lambda^t$. By definition of $W_h\backslash\big(S^\delta W_h\big)$, we have 
\[
W_h\big(W_h\backslash S^\delta W_h\big)\leq S^\delta W_h\enspace ,
\]
hence
\begin{eqnarray*}
W_h\mu^{t+\delta} & = & W_h\big(W_h\backslash S^\delta W_h\big)\mu^t\\
& \leq & \big(S^\delta W_h\big)\mu^t\\
& \leq & S^\delta W_h\lambda^t\enspace .
\end{eqnarray*}
Since 
\begin{eqnarray*}
\lambda^{t+\delta} & = & W_h\backslash\big(S^\delta W_h\lambda^t\big)\\
& = & \max\set{\lambda\in\rmaxb^p}{W_h\lambda\leq S^\delta W_h\lambda^t}\enspace ,
\end{eqnarray*} 
we get that $\mu^{t+\delta}\leq\lambda^{t+\delta}$. 
Then $\mu^{t}\leq\lambda^{t}$ for $t=0,\delta,...,T$. Since $W_h$ is monotone, we deduce~\eqref{comparaison}.
\end{proof}
An approximation of \eqref{mut+dt} using formulae of the same type as
\eqref{stilde} is also discussed in~\cite{mceneaney-hortonb}.
\section{Error analysis}\label{section_error}
\subsection{General error estimates}
In the sequel we denote by $\|v\|_{\infty}=\sup_{i\in I}|v(i)|\in\R\cup\{+\infty\}$ the sup-norm of any function $v:I\to\R$. We also use the same notation $\|v\|_{\infty}=\sup_{i\in I}|v_i|$ for a vector $v=(v_i)_{i\in I}$. For any two sets $I$ and $J$, a map $\Phi:\R^I\to\R^J$ is said monotone and homogeneous if it is monotone for the natural order and if for all $u\in\R^I$ and $\lambda\in\R$, $\Phi(u+\lambda)=\Phi(u)+\lambda$ with $(u+\lambda)(i)=u(i)+\lambda$. Monotone homogeneous maps are nonexpansive for the sup-norm: $\|\Phi(u)-\Phi(v)\|_{\infty}\leq\|u-v\|_{\infty}$, see~\cite{crandall}. In particular, max-plus or min-plus linear operators are non-expansive for the sup-norm. This property will be frequently used in the sequel. In order to simplify notations, we denote $\bar\tau_{\delta}=\{0,\delta,\cdots,T\}$, $\tau_{\delta}^+=\bar\tau_\delta\backslash\{0\}$ and $\tau_{\delta}^-=\bar\tau_\delta\backslash\{T\}$ .
\begin{rem}
To establish the main result of the paper (Theorem~\ref{th-main} below), we shall need only to take the norm of finite valued functions. However, we wish to emphasize that all the computations that follow are valid for functions with values in $\rbar$ if one replaces every occurence of a term of the form $\|u-v\|_\infty$ by $d_\infty(u,v)=\inf\set{\lambda\geq 0}{-\lambda+v\leq u\leq\lambda+v}$. Observe that $d_\infty(u,v)$ is a semidistance and that $d_\infty(u,v)=\|u-v\|_\infty$, if $u-v$ takes finite values. Observe also that if a map $\Phi:\rbar^I\to\rbar^J$ is monotone and homogeneous, $d_\infty(\Phi(u),\Phi(v))\leq d_\infty(u,v)$, for all $u,v\in\rbar^I$.
\end{rem}
The following lemma shows that the error of the ideal max-plus finite element method is controlled by the projection errors $\|\projimker{W_h}{Z_h^*}(v^t)-v^t\|_{\infty}$. This lemma may be thought of as an analogue of Cea's lemma in the classical analysis of the errors of the finite element method. Projectors over semimodules in the MFEM correspond to orthogonal projectors in the classical finite element method. 
\begin{lem}\label{lemme-general}
For $t\in\bar\tau_{\delta}$, let $v^t$ be the value function at time $t$, and $v_h^t$ be its approximation given by the ideal max-plus finite element method. We have
\begin{equation}
\|v_h^T-v^T\|_{\infty}\leq \|\projimker{W_h}{}(v^0)-v^0\|_{\infty}+\sum_{t\in\tau_{\delta}^+}\|\projimker{W_h}{Z_h^*}(v^t)-v^t\|_{\infty}\enspace .
\end{equation}
\end{lem}
\begin{proof}
For all $t\in\tau_\delta^-$, we have
\begin{eqnarray*}
\|v_h^{t+\delta}-v^{t+\delta}\|_{\infty}
&\leq& \|v_h^{t+\delta}-S_h^\delta (v^t)\|_{\infty}+\|S_h^\delta (v^t)-v^{t+\delta}\|_{\infty}\\
&\leq& \|S_h^\delta (v_h^t)-S_h^\delta (v^t)\|_{\infty}+\|\projimker{W_h}{Z_h^*}\comp S^\delta (v^t)-v^{t+\delta}\|_{\infty}\enspace .
\end{eqnarray*}
Since $S_h^\delta$ is a non-expansive operator, we deduce
\begin{equation*}
\|v_h^{t+\delta}-v^{t+\delta}\|_{\infty}\leq\|v_h^t-v^t\|_{\infty}+\|\projimker{W_h}{Z_h^*}(v^{t+\delta})-v^{t+\delta}\|_{\infty}\enspace .
\end{equation*}
The result is obtained by induction on $t$, using the fact that $v_h^0=P_{\calw_h}(v^0)=\projimker{W_h}{}(v^0)$.
\end{proof}
To obtain an error estimate, we need to bound $\|\projimker{W_h}{Z_h^*}(v^{t})-v^{t}\|_{\infty}$ for all $t\in\tau_{\delta}^+$. Since $\projimker{W_h}{Z_h^*}=\projimker{W_h}{}\comp\projimker{}{Z_h^*}$, we have
\begin{eqnarray*}
\|\projimker{W_h}{Z_h^*}(v^t)-v^t\|_{\infty} &=& \|\projimker{W_h}{}\comp\projimker{}{Z_h^*}(v^t)-v^t\|_{\infty}\\
&\leq& \|\projimker{W_h}{}\comp\projimker{}{Z_h^*}(v^t)-\projimker{W_h}{}(v^t)\|_{\infty}+\|\projimker{W_h}{}(v^t)-v^t\|_{\infty}\enspace ,
\end{eqnarray*}
and since $\projimker{W_h}{}$ is a non-expansive operator, we get
\begin{equation}
\|\projimker{W_h}{Z_h^*}(v^t)-v^t\|_{\infty}\leq\|\projimker{}{Z_h^*}(v^t)-v^t\|_{\infty}+\|\projimker{W_h}{}(v^t)-v^t\|_{\infty}\enspace .
\end{equation}
Using this inequality together with Lemma~\ref{lemme-general}, we deduce the following corollary.
\begin{cor}\label{err-proj}
For $t\in\bar\tau_{\delta}$, let $v^t$ be the value function at time $t$, and $v_h^t$ be its approximation given by the ideal max-plus finite element method. We have
\begin{equation*}
\|v_h^T-v^T\|_{\infty}\leq (1+\frac T \delta)\Big(\sup_{t\in\bar\tau_\delta}(\|\projimker{}{Z_h^*}(v^t)-v^t\|_{\infty}+\|\projimker{W_h}{}(v^t)-v^t\|_{\infty})\Big)\enspace .
\end{equation*}
\end{cor}
The following general lemma shows that the error of the effective finite element
method is controlled by the projection errors and the errors resulting from the approximation of the matrix $K_h$ by a matrix $\tilde K_h$.
\begin{lem}\label{error_effect1}
For $t\in\bar\tau_\delta$,
let $v^t$ be the value function at time $t$, and $v_h^t$
be its approximation given by the effective max-plus finite element method, where $K_h$ is approximated by $\tilde K_h$. We have
\begin{align*}
\|v_h^T- v^T\|_{\infty}&\leq (1+\frac{T}{\delta})\Big(\sup_{t\in\bar\tau_\delta}(\|\projimker{}{Z_h^*}(v^t)-v^t\|_{\infty}+\|\projimker{W_h}{}(v^t)-v^t\|_{\infty})\\
&\qquad\qquad\qquad +\|\tilde K_h-K_h\|_{\infty}\Big)\enspace .
\end{align*}
\end{lem}
\begin{proof}
Since $v_h^t$ is computed with the approximation $\tilde K_h$ of $K_h$, we have $v_h^t=W_h\lambda^t$, $t\in\bar\tau_\delta$, with 
\[
\lambda^{t+\delta}=M_h\sh\comp(\tilde K_h\lambda^t)=W_h\sh\comp (Z_h^*)\sh\comp(\tilde K_h\lambda^t)\enspace .
\]
We have
\begin{eqnarray*}
\|v_h^{t+\delta}-v^{t+\delta}\|_{\infty} &\leq& \|v_h^{t+\delta}-S_h^\delta v_h^t\|_{\infty}+\|S_h^\delta v_h^t-S_h^\delta v^t\|_{\infty}+\|S_h^\delta v^t-v^{t+\delta}\|_{\infty}\\
&\leq& \|\projimker{W_h}{}\comp(Z_h^*)\sh\comp(\tilde K_h\lambda^t)-\projimker{W_h}{}\comp(Z_h^*)\sh\comp Z_h^*\comp S^\delta W_h\lambda^t\|_{\infty}\\
&& \qquad +\|v_h^t-v^t\|_{\infty}+\|\projimker{W_h}{Z_h^*}(v^{t+\delta})-v^{t+\delta}\|_{\infty}\\
&\leq& \|\tilde K_h\lambda^t-K_h\lambda^t\|_{\infty}+\|v_h^t-v^t\|_{\infty}+\|\projimker{W_h}{Z_h^*}(v^{t+\delta})-v^{t+\delta}\|_{\infty}\\
&\leq& \max_{\begin{subarray}{c}
1\leq j\leq q\\
1\leq i\leq p
\end{subarray}}|(\tilde K_h)_{ji}-(K_h)_{ji}|
+\|v_h^t-v^t\|_{\infty}+\|\projimker{W_h}{Z_h^*}(v^{t+\delta})-v^{t+\delta}\|_{\infty}\enspace .
\end{eqnarray*}
We deduce that
\begin{equation*}
\|v_h^T-v^T\|_{\infty}\leq\|\projimker{W_h}{}(v^0)-v^0\|_{\infty}+\sum_{t\in\tau_\delta^+}\Big(\|\projimker{W_h}{Z_h^*}(v^t)-v^t\|_{\infty}+\|\tilde K_h-K_h\|_{\infty}\Big)\enspace ,
\end{equation*}
and so
\begin{align*}
\|v_h^T- v^T\|_{\infty}&\leq (1+\frac{T}{\delta})\Big(\sup_{t\in\bar\tau_\delta}(\|\projimker{}{Z_h^*}(v^t)-v^t\|_{\infty}+\|\projimker{W_h}{}(v^t)-v^t\|_{\infty})\\
&\qquad\qquad\qquad+\|\tilde K_h-K_h\|_{\infty}\Big)\enspace .
\end{align*}
\end{proof}
\begin{cor}\label{error_effect}
For $t\in\bar\tau_\delta$,
let $v^t$ be the value function at time $t$, and $v_h^t$
be its approximation given by the effective max-plus finite element method,
implemented with the approximation $K_{H,h}$ of $K_h$,
given by~\eqref{appr_K}.
We have
\begin{align*}
\|v_h^T- v^T\|_{\infty}&\leq (1+\frac{T}{\delta})\Big(\sup_{t\in\bar\tau_\delta}(\|\projimker{}{Z_h^*}(v^t)-v^t\|_{\infty}+\|\projimker{W_h}{}(v^t)-v^t\|_{\infty})\\
&\qquad\qquad\qquad+\max_{1\leq i\leq p}\|[S^\delta w_i]_H-S^\delta w_i\|_{\infty}\Big)\enspace .
\end{align*}
\end{cor}
\begin{proof}
Using the same technique as in the precedent lemma and using that $K_{H,h}=Z_h^*[S^\delta W_h]_H$ and $K_h=Z_h^*S^\delta W_h$ we have
\begin{align}\label{error-K}
\|K_{H,h}-K_h\|_{\infty}&\leq\|[S^\delta W_h]_H-S^\delta W_h\|_{\infty}\nonumber\\
&=\max_{1\leq i\leq p}\|[S^{\delta}w_i]_H -S^{\delta}w_i\|_{\infty}\enspace ,
\end{align}
which ends the proof.
\end{proof}
\begin{cor}\label{complete-error}
For $t\in\bar\tau_\delta$,
let $v^t$ be the value function at time $t$, and $v_h^t$
be its approximation given by the effective max-plus finite element method,
implemented with the approximation $\tilde K_{H,h}$ of $K_h$,
given by~\eqref{e-convenient}.
We have
\begin{align*}
\|v_h^T- v^T\|_{\infty}&\leq (1+\frac{T}{\delta})\Big(\sup_{t\in\bar\tau_\delta}(\|\projimker{}{Z_h^*}v^t-v^t\|_{\infty}+\|\projimker{W_h}{}v^t-v^t\|_{\infty})\\
&\qquad \qquad\qquad+\max_{1\leq i\leq p}\|[S^{\delta}w_i]_H -S^{\delta}w_i\|_{\infty}+\|\tilde K_{H,h}-K_{H,h}\|_{\infty}\Big)\enspace .
\end{align*}
\end{cor}
\begin{proof}
We use Lemma~\ref{error_effect1}, together with Equation~\eqref{error-K} and
\[
\|\tilde K_{H,h}-K_h\|_{\infty}\leq \|\tilde K_{H,h}-K_{H,h}\|_{\infty}+\|K_{H,h}-K_h\|_{\infty}\enspace .
\]
\end{proof}
\subsection{Projection errors}
In this section, we estimate the projection errors resulting from different choices of finite elements. 
Recall that a function $f$ is $c$\new{-semiconvex} if $f(x)+\frac c 2\|x\|_2^2$, where $\|\cdot\|_2$ is the standard euclidean norm of $\R^n$, is convex. A function $f$ is $c$\new{-semiconcave} if $-f$ is $c$-semiconvex. Spaces of semiconvex functions were intensively used in the max-plus based approximation method of Fleming and McEneaney~\cite{mceneaney}, see also~\cite{mceneaney-hortona}, \cite{mceneaney-hortonb}, \cite{mceneaney02}, \cite{mceneaney03}, \cite{mceneaney04}, \cite{falcone}, \cite{capuzzodolcetta-ishii}, \cite{capuzzodolcetta-falcone}.

We shall use the following finite elements.
\begin{defin}[$P_1$ finite elements]
We call \new{$P_1$ finite element} or \new{Lipschitz finite element} centered at point $\hat x\in X$,
with constant $a>0$, the function $w(x)=-a\|x-\hat x\|_1$ where $\|x\|_1=\sum_{i=1}^n |x_i|$ is the $l^1$-norm of $\R^n$.
\end{defin}
The family of Lipschitz finite element of constant $a$ generates, in the max-plus sense, the semimodule of Lipschitz continuous functions from $X$ to $\bar\R$ of Lipschitz constant $a$ with respect to $\|\cdot\|_1$.
\begin{defin}[$P_2$ finite elements]
We call \new{$P_2$ finite element} or \new{quadratic finite element} centered at point $\hat x\in X$,
with Hessian $c>0$, the function $w(x)=-\frac{c}{2}\|x-\hat x\|_2^2$.
\end{defin}
When $X=\R^n$, the family of quadratic finite elements with Hessian $c$ generates, in the max-plus sense, the semi-module of lower-semicontinuous $c$-semiconvex functions with values in $\bar\R$.\\
\textbf{Notations.} Let $Y$ be a subset of $\R^n$ and $f$ be a function from $Y$ to $\rbar$. We will denote by $\conv Y$ the convex hull of $Y$, $\ri Y$ the relative interior of $Y$, $\dom f$ the effective domain of $f$ and $\partial f(x)$ the subdifferential of $f$ at $x\in\dom f$. \\
When $C$ is a nonempty convex subset of $\R^n$ and $c>0$, a fonction is said to be $c$-\new{strongly convex} on $C$ if and only if $f-\frac{1}{2}c\|\cdot\|_2^2$ is convex on $C$. A function $f$ is $c$-\new{strongly concave} on $C$ if $-f$ is $c$-strongly convex on $C$.

\begin{figure}[htbp]
\centering
\input{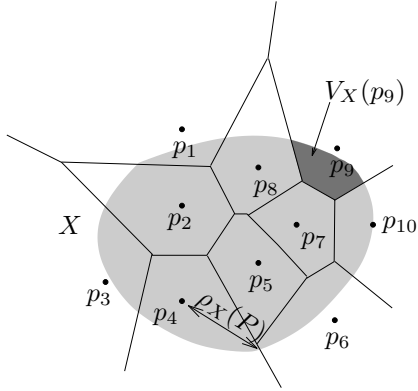}
\caption{Voronoi tesselation}
\label{voronoi}
\end{figure}

Let $P$ be a finite subset of $\R^n$. The \new{Voronoi cell} of a point $p\in P$ is defined by
\[
V(p)=\set{x\in \R^n}{\|x-p\|_2\leq\|x-q\|_2,\forall q\in P}.
\]
The family $\{V(p)\}_{p\in P}$ constitutes a subdivison of $\R^n$, which is called a Voronoi tesselation (see~\cite{comp-geom} for an introduction to Voronoi tesselations). We define \new{the restriction of $V(p)$ to $X$} to be:
\[
V_X(p)=V(p)\cap X.
\]
We define $\rho_X(P)$ to be the \new{maximal radius} of the restriction to X of the Voronoi cells of the points of $P$:
\[
\rho_X(P):=\sup_{p\in P}\sup_{x\in V_X(p)}\|x-p\|_2.
\]
Observe that 
\[
\rho_X(P):=\sup_{x\in X}\inf_{p\in P}\|x-p\|_2.
\]
The previous definitions are illustrated in Figure~\ref{voronoi}. The set $X$ is in light gray, $P=\{p_1,\cdots, p_{10}\}$, $V_X(P_9)$ is in dark gray and $\rho_X(P)$ is indicated by a bidirectional arrow.

The next two lemmas bound the projection error in term of the radius of Voronoi cells.
\begin{lem}[Primal projection error]\label{proj-p2}
Let $X$ be a compact convex subset of $\R^n$. Let $v:X\to \R$ be $c$-semiconvex and Lipschitz continuous function with Lipschitz constant $L_v$ with respect to the euclidean norm. Let $v_c(x)=v(x)+\frac c 2\|x\|_2^2$. Let $\hat X=X+\mathrm{B}_2(0,\frac{L_v}{c})$, let $\hat X_h$ be a finite subset of $\R^n$ and let $\calw_h$ denote the complete subsemimodule of $\rmaxb^{_{\scriptstyle X}}$ generated by the family $(w_{\hat x_h})_{\hat x_h\in\hat X_h}$ where $w_{\hat x_h}(x)=-\frac{c}{2}\|x-\hat x_h\|_2^2$. Then 
\[
\|v-P_{\calw_h} v\|_{\infty}\leq c\operatorname{diam}X\rho_{\hat X}(\hat X_h)
\]
\end{lem}
\begin{proof}
Let $\calw$ denote the complete subsemimodule of $\rmaxb^{_{\scriptstyle X}}$ generated by the family $(w_{\hat x})_{\hat x\in\hat X}$. We will first prove that for all $x\in X$, $P_\calw v(x)=v(x)$. It is obvious that $\forall x\in X$, $v(x)\geq P_\calw v(x)$. Using that $P_{\calw}=W\comp W\sh$, with $W=\mathrm{col}(w_{\hat x})_{\hat x\in \hat X}$, we obtain
\begin{eqnarray*}
P_\calw v(x) &=& \sup_{\hat x\in\hat X}\Big(-\frac{c}{2}\|x-\hat x\|_2^2+\inf_{y\in X}\big(\frac{c}{2}\|y-\hat x\|_2^2+v(y)\big)\Big)\\
&=& \sup_{\hat x\in\hat X}\Big(-\frac{c}{2}\|x\|_2^2+c\hat x\cdot x-\sup_{y\in X}\big(c\hat x\cdot y-v_c(y)\big)\Big)\\
&=& -\frac{c}{2}\|x\|_2^2+\sup_{\hat x\in\hat X}\big(c\hat x\cdot x-v_c^*(c\hat x)\big)\enspace ,
\end{eqnarray*}
where $v_c^*$ denotes the Fenchel transform of $v_c$. Since $v_c$ is l.s.c., convex and proper, we have for all $x\in X$
\begin{equation}\label{fenchel}
v_c(x)=v_c^{**}(x)=\sup_{\theta\in\R^n}\big(\theta\cdot x-v_c^*(\theta)\big)\enspace .
\end{equation} 
Using Theorem $23.4$ of~\cite{rockafellar}, for all $x\in\ri(\dom v_c)$, the subdifferential of $v_c$ at $x$, $\partial v_c(x)=\set{\theta\in\R^n}{v_c(y)-v_c(x)\geq\theta\cdot(y-x),\forall y\in X}$, is non-empty. Then $\theta\in\partial v_c(x)$ if and only if $v_c^*(\theta)=\theta\cdot x-v_c(x)$ and consequently, the supremum of~\eqref{fenchel} is attained for all elements $\theta$ of $\partial v_c(x)$.\\
Set $q(x)=\frac{c}{2}\|x\|_2^2$. Using the fact that $q(y)-q(x)=q'(x)\cdot(y-x)+O(\|y-x\|_2^2)$ and that $v$ is Lipschitz continuous with Lipschitz constant $L_v$, we obtain $\partial v_c(x)\subset\mathrm{B}_2(cx,L_v)$ for all $x\in\ri X$. Therefore, for all $x\in\ri X$,
\begin{equation}\label{v_c}
v_c(x)=\sup_{\hat x\in\hat X}\big(c\hat x\cdot x-v_c^*(c\hat x)\big)\enspace .
\end{equation}
By continuity in the members of Equation~\eqref{v_c}, we have the equality for all $x\in X$, and so 
\begin{eqnarray*}
P_\calw v(x) &=& -\frac{c}{2}\|x\|_2^2+\sup_{\hat x\in\hat X}\big(c\hat x\cdot x-v_c^*(c\hat x)\big)\\
&=& -\frac{c}{2}\|x\|_2^2+v_c(x)\\
&=& v(x)\enspace ,
\end{eqnarray*}
for all $x\in X$.\\
Now, fix $x\in X$. For $\hat x\in\hat X$, we set $\varphi(\hat x)=c\hat x\cdot x-v_c^*(c\hat x)$. Since $P_{\calw_h}v\leq P_\calw v\leq v$, we have for all $x\in X$
\begin{eqnarray*}
0\leq v(x)-P_{\calw_h}v(x) &=& P_\calw v(x)-P_{\calw_h}v(x)\\
&=& \sup_{\hat x\in\hat X}\varphi(\hat x)-\sup_{\hat x_h\in\hat X_h}\varphi(\hat x_h)\\
&=&\sup_{\hat x\in\hat X}\inf_{\hat x_h\in\hat X_h}\varphi(\hat x)-\varphi(\hat x_h)\enspace .
\end{eqnarray*}
We have $\partial(-\varphi)(\hat x)=-cx+c\partial v_c^*(c\hat x)$. Since $\partial v_c^*\subset X$, we have $\partial(-\varphi)(\hat x)\subset c(X-x)\subset\mathrm{B}_2(0,c\operatorname{diam}X)$. Hence, $\varphi$ is Lipschitz continuous with Lipschitz constant $L_\varphi=c\operatorname{diam}X$. Then for all $x\in X$
\begin{eqnarray*}
v(x)-P_{\calw_h}v(x) &\leq& \sup_{\hat x\in\hat X}\inf_{\hat x_h\in\hat X_h}L_\varphi\|\hat x-\hat x_h\|_2\\
&=& c\operatorname{diam}X\rho_{\hat X}(\hat X_h)\enspace .
\end{eqnarray*}
\end{proof}
\begin{lem}[Dual projection error]\label{proj-p1}
Let $X$ be a bounded subset of $\R^n$ and $\hat X$ a finite subset of $\R^n$. Let $v:X\to\R$ be a given Lipschitz continuous function with Lipschitz constant $L_v$ with respect to the euclidean norm. Let $\calz_{\hat X}$ denote the complete semimodule of $\rmaxb^{_X}$ generated by the $P_1$ finite elements $(z_{\hat x})_{\hat x\in\hat X}$ centered at the points of $\hat X$ with constant $a\geq L_v$. Then 
\[
\|v-P^{-(\calz_{\hat X})}v\|_\infty\leq n(a+L_v)\rho_X(\hat X).
\]
\end{lem}
\begin{proof}
It is clear that $P^{-(\calz_{\hat X})}v\geq v$ and using that $P^{-(\calz_{\hat X})}=(Z^*)\sh\comp Z^*$, with $Z=\mathrm{col}(z_{\hat x})_{\hat x\in\hat X}$, we obtain
\begin{equation*}
P^{-(\calz_{\hat X})}v(x)-v(x) = \inf_{\hat x\in\hat X}\Big(a\|x-\hat x\|_1+\sup_{y\in X}\big(-a\|y-\hat x\|_1+v(y)-v(x)\big)\Big)\enspace ,
\end{equation*}
for all $x\in X$. Since $v$ is $L_v$-Lipschitz continuous, we have
\begin{eqnarray*}
P^{-(\calz_{\hat X})}v(x)-v(x) &\leq& \inf_{\hat x\in\hat X}\Big(a\|x-\hat x\|_1+\sup_{y\in X}\big(-a\|y-\hat x\|_1+L_v\|y-x\|_2\big)\Big)\\
&\leq& \inf_{\hat x\in\hat X}\Big(a\|x-\hat x\|_1+\sup_{y\in X}\big(-a\|y-\hat x\|_1+L_v\|y-x\|_1\big)\Big)\\
&\leq& \inf_{\hat x\in\hat X}\Big(a\|x-\hat x\|_1+\sup_{y\in X}\big(-a\|y-\hat x\|_1+L_v\|y-\hat x\|_1\\
&&\qquad\qquad +L_v\|x-\hat x\|_1\big)\Big)\\
&=& \inf_{\hat x\in\hat X}\Big((a+L_v)\|x-\hat x\|_1+\sup_{y\in X}(L_v-a)\|y-\hat x\|_1\Big)\enspace .
\end{eqnarray*}
Since $a\geq L_v$, we deduce
\begin{equation*}
P^{-(\calz_{\hat X})}v(x)-v(x) \leq (a+L_v)\sup_{x\in X}\inf_{\hat x\in\hat X}\|x-\hat x\|_1\leq n(a+L_v)\rho_X(\hat X)\enspace .
\end{equation*}
\end{proof}
\subsection{The approximation errors}
To state an error estimate, we make the following standard assumptions (see~\cite{barles} for instance):
\begin{itemize}
\item[-] $(H1)$ $f:X\times U\to \R^n$ is bounded and Lipschitz continuous with
  respect to $x$, meaning that there exist $L_f >0$ and $M_f >0$ such that
\begin{align*}
\| f(x,u)-f(y,u)\|_2 &\leq L_f\| x-y\|_2   &\forall x,y\in X,u \in U,\\
\| f(x,u)\|_2 &\leq M_f, &\forall x\in X, u\in U\enspace .
\end{align*}
\item[-] $(H2)$ $\ell:X\times U\to \R$ is bounded and Lipschitz continuous with respect to $x$, meaning that there exist $L_\ell>0$ and $M_\ell>0$ such that
\begin{align*}
| \ell(x,u)-\ell(y,u) |&\leq L_\ell\| x-y\|_2 & \forall x,y \in X, u \in U,\\
| \ell(x,u) |&\leq M_\ell, & \forall  x\in X, u \in U\enspace .
\end{align*}
\end{itemize}
\subsubsection{Approximation of $S^\delta w$}
\begin{lem}\label{appr-StI}
Let $X$ be a convex subset of $\R^n$. We
make assumptions \mrm{(H1)} and \mrm{(H2)}. Let $w:x\to\R$ be such that $w$ is $\mathcal{C}^1$ on a neighborhood of $X$, Lipschitz continuous with Lipschitz constant $L_w$ with respect to the euclidean norm, $c_1$-semiconvex and $c_2$-semiconcave. Then there exists $K_1>0$ such that $\|[S^\delta w]_H-S^\delta w\|_{\infty}\leq K_1\delta^2$, for $\delta>0$, where $[S^\delta w]_H$ is given by~\eqref{stilde}.
\end{lem}
\begin{proof}
We first show that there exists $K_1>0$ such that 
\[
[S^\delta w]_H(x)-S^\delta w(x)\geq -K_1\delta^2, \quad \forall x\in X\enspace .
\]
For all $x\in X$ and $u\in U$, define $\mb x_{u,x}$ to be the trajectory
such that $\dot{\mb x}_{u,x}(s)=f(\mb x_{u,x}(s),u)$ and $\mb x_{u,x}(0)=x$.
In other words, we apply a constant control $u$.
We have
\begin{equation*}
(S^{\delta}w)(x)\geq\sup\set{\int_0^{\delta}\ell(\mb x_{u,x}(s),u)ds+w(\mb x_{u,x}(\delta))}{u\in U} \enspace .
\enspace.
\end{equation*}
Since $\ell$ is Lipschitz continuous and $f$ is bounded, we have
\begin{eqnarray*}
\Big|\int_0^{\delta}[\ell(\mb x_{u,x}(s),u)-\ell(x,u)]ds\Big| &\leq& L_\ell\int_0^\delta\|\mb x_{u,x}(s)-x\|_2ds\\ 
&\leq& L_\ell\int_0^\delta M_f sds\enspace ,
\end{eqnarray*}
then
\begin{equation}\label{l_lip}
\Big|\int_0^{\delta}[\ell(\mb x_{u,x}(s),u)-\ell(x,u)]ds\Big|\leq\frac{1}{2}L_{\ell}M_f\delta^2\enspace .
\end{equation}
Therefore
\begin{align*}
(S^{\delta}w)(x)&\geq-\frac{1}{2}L_{\ell}M_f\delta^2+\sup\set{\delta\ell(x,u)+w(\mb x_{u,x}(\delta))}{u\in U}\enspace .
\end{align*}
Since $w$ is Lipschitz continuous and $f$ is bounded and Lipschitz continuous, we have
\begin{eqnarray*}
\Big|w(\mb x_{u,x}(\delta))-w(x+\delta f(x,u))\Big|&\leq& L_w\|\mb x_{u,x}(\delta)-x-\delta f(x,u)\|_2\\ 
&\leq& L_w\int_0^\delta\|f(\mb x_{u,x}(s),u)-f(x,u)\|_2 ds\\
&\leq& L_w\int_0^\delta L_f\|\mb x_{u,x}(s)-x\|_2ds\\
&\leq& L_wL_f\int_0^\delta M_fsds\enspace ,
\end{eqnarray*}
and so
\begin{equation}\label{w_lip}
\Big|w(\mb x_{u,x}(\delta))-w(x+\delta f(x,u))\Big|\leq\frac 1 2 L_wL_fM_f\delta^2\enspace .
\end{equation}
Moreover, since $w$ is $c_1$-semiconvex, we have
\begin{equation}\label{w_c1semiconvex}
w(x+\delta f(x,u))\geq w(x)+\delta\nabla w(x)\cdot f(x,u)-\frac{c_1}{2}M_f^2\delta^2\enspace .
\end{equation}
We deduce from~\eqref{l_lip}, \eqref{w_lip} and~\eqref{w_c1semiconvex} 
\begin{eqnarray*}
(S^{\delta}w)(x)&\geq&-\big(L_{\ell}M_f+L_wL_fM_f+c_1M_f^2\big)\frac{\delta^2}{2}+w(x)\\
&&\qquad\qquad\qquad+\sup_{u\in U}\big\{\delta\ell(x,u)+\delta\nabla w(x)\cdot f(x,u)\big\}\\
&\geq& -\big(L_{\ell}M_f+L_wL_fM_f+c_1M_f^2\big)\frac{\delta^2}{2}+w(x)+\delta H(x,\nabla w(x))\enspace .
\end{eqnarray*}
This ends the first part of the proof.

We now prove an opposite inequality. 
For all $x\in X$ and for all measurable
functions $\mb u:[0,\delta]\to U$, define
$\mb x_{\mb u,x}$ to be the trajectory
such that $\dot{\mb x}_{\mb u,x}(s)=f(\mb x_{\mb u,x}(s),\mb u(s))$ and 
$\mb x_{\mb u,x}(0)=x$.
Since $\ell(x,u)\leq H(x,p)-p\cdot f(x,u)$, 
for all $p\in\R^n$, $x\in X$ and $u\in U$, we deduce that
\begin{align*}
(S^\delta w)(x)&\leq
\sup\Big\{\int_0^\delta H(\mb x_{\mb u,x}(s),\nabla w(x))ds+
w(\mb x_{\mb u,x}(\delta))\\
&\qquad \qquad \qquad\qquad-\nabla w(x)\cdot\int_0^\delta f(\mb x_{\mb u,x}(s),\mb u(s))ds\mid 
 \mb u: [0,\delta] \to U\Big\}\\
&=\sup\Big\{\int_0^\delta H(\mb x_{\mb u,x}(s),\nabla w(x))ds\\
&\qquad \qquad \qquad \qquad+w(\mb x(\delta))-\nabla w(x)\cdot\big(\mb x_{\mb u,x}(\delta)-x\big)\mid\mb u: [0,\delta] \to U\Big\}\enspace .
\end{align*}
Using the fact that $\ell$ and $f$ are Lipschitz continuous with respect to $x$, we have for all $x,x'\in X$, $p\in\R^n$
\begin{equation*}
\Big|H(x,p)-H(x',p)\Big|\leq \big(L_\ell+L_f\|p\|_2\big)\|x-x'\|_2\enspace ,
\end{equation*}
therefore
\begin{eqnarray*}
(S^\delta w)(x)&\leq&\sup\Big\{\big(L_\ell+L_fL_w\big)\int_0^\delta\|\mb x_{\mb u,x}(s)-x\|_2ds+\delta H(x,\nabla w(x))\\
&&\qquad\qquad+w(\mb x_{\mb u,x}(\delta))-\nabla w(x)\cdot\big(\mb x_{\mb u,x}(\delta)-x\big)\mid\mb u:[0,\delta]\to U\Big\}\\
&\leq&\big(L_\ell+L_fL_w\big)M_f\frac{\delta^2}{2}+\delta H(x,\nabla w(x))\\
&&\qquad+\sup\Big\{w(\mb x_{\mb u,x}(\delta))-\nabla w(x)\cdot\big(\mb x_{\mb u,x}(\delta)-x\big)\mid\mb u:[0,\delta]\to U\Big\}\enspace .
\end{eqnarray*}
Since $w$ is $c_2$-semiconcave, we have
\begin{equation*}
w(\mb x_{\mb u,x}(\delta))\leq w(x)+\nabla w(x)\cdot\big(\mb x_{\mb u,x}(\delta)-x\big)+\frac{c_2}{2}M_f^2\delta^2\enspace .
\end{equation*}
We obtain
\begin{equation*}
(S^\delta w)(x)\leq \big(L_\ell+L_fM_{Dw}+c_2M_f\big)M_f\frac{\delta^2}{2}+w(x)+\delta H(x,\nabla w(x))\enspace .
\end{equation*}
To end the proof, we take $K_1=\frac{1}{2}\big(L_\ell M_f+L_fM_{Dw}M_f+\max(c_1,c_2)M_f^2\big)$.
\end{proof}
\subsubsection{Approximation of the matrix $K_h$ by the matrix $\tilde K_H$}
\begin{lem}
Let $X$ be a compact subset of $\R^n$. We consider  an upper semicontinuous function $\varphi:X\to\R$ and a Lipschitz continuous function $\psi:X\to\R$ with Lipschitz constant $L_\psi$ with respect to a norm $\|\cdot\|$. For $\varepsilon\geq 0$, we define:
\begin{subequations}\label{fepsilon}
\begin{gather}
F_{\varepsilon}=\set{x\in X}{\varphi(x)\geq\sup_{x'\in X}\varphi(x')-\varepsilon},\\
g(\varepsilon)=\sup_{x\in F_{\varepsilon}}d(x,F_0)\enspace ,
\end{gather}
\end{subequations}
where $d(x,F_0)=\inf_{y\in F_0}\|y-x\|$. We have:
\begin{equation*}
\mid\sup_{x\in X}\big(\varphi(x)+\delta\psi(x)\big)-\big[\sup_{x\in
    X}\varphi(x)+\delta\sup_{x\in
    \argmax{\varphi}}\psi(x)\big]\mid\leq L_\psi\delta g(\delta M)\enspace ,
\end{equation*}
where $M=\sup_{x\in X}\psi(x)-\inf_{x\in X}\psi(x)$.
\end{lem}
\begin{proof}
Since $\varphi$ is u.s.c.\ and $X$ is compact, $F_0=\argmax\varphi$ and
\begin{equation}\label{+grand0}
\sup_{x\in X}\big(\varphi(x)+\delta\psi(x)\big)\geq \sup_{x\in X}\varphi(x)+\delta\sup_{x\in F_0}\psi(x)\enspace .
\end{equation}
For $\varepsilon>0$ we have:
\begin{equation*}
\sup_{x\in X}\big(\varphi(x)+\delta\psi(x)\big)=\max\big[\sup_{x\in F_{\varepsilon}}\big(\varphi(x)+\delta\psi(x)\big), \sup_{x\in X\backslash F_{\varepsilon}}\big(\varphi(x)+\delta\psi(x)\big)\big]\enspace .
\end{equation*}
Let $\varepsilon=\delta (\sup_{x\in X}\psi(x)-\inf_{x\in X}\psi(x))=M\delta$ (which is finite since $\psi$ is continuous and $X$ is compact). We have:
\begin{eqnarray*}
\sup_{x\in X\backslash F_{\varepsilon}}\big(\varphi(x)+\delta\psi(x)\big)&\leq&-\varepsilon+\sup_{x\in X}\varphi(x)+\delta\sup_{x\in X}\psi(x)\\
&=&\sup_{x\in F_{\varepsilon}}\varphi(x)+\delta\inf_{x\in X}\psi(x)\\
&\leq&\sup_{x\in F_{\varepsilon}}\big[\varphi(x)+\delta\psi(x)\big]\enspace .
\end{eqnarray*}
Therefore
\begin{eqnarray}\label{+petit}
\sup_{x\in X}\big(\varphi(x)+\delta\psi(x)\big)&=&\sup_{x\in
  F_{\varepsilon}}\big(\varphi(x)+\delta\psi(x)\big)\nonumber\\
&\leq&\sup_{x\in X}\varphi(x)+\delta\sup_{x\in F_{\varepsilon}}\psi(x)\enspace .
\end{eqnarray}
We deduce from \eqref{+grand0} and \eqref{+petit}:
\begin{equation*}
0\leq\sup_{x\in X}\big(\varphi(x)+\delta\psi(x)\big)-\big[\sup_{x\in
    X}\varphi(x)+\delta\sup_{x\in F_0}\psi(x)\big]\leq\delta\big[\sup_{x\in F_{\varepsilon}}\psi(x)-\sup_{x\in F_0}\psi(x)\big]\enspace .
\end{equation*}
Since $\psi$ is Lipschitz continuous, we have
\begin{eqnarray*}
 \sup_{x\in F_{\varepsilon}}\psi(x)-\sup_{x\in F_0}\psi(x)&=&\sup_{x\in F_{\varepsilon}}\inf_{y\in F_0}\big(\psi(x)-\psi(y)\big)\\
&\leq&\sup_{x\in F_{\varepsilon}}\inf_{y\in F_0}L_\psi\|x-y\|\\
&=&L_\psi\sup_{x\in F_{\varepsilon}}d(x,F_0)\\
&=&L_\psi g(\varepsilon)\enspace .
\end{eqnarray*}
\end{proof}
\begin{cor}\label{erreurB}
Let $X$ be a compact convex subset of $\R^n$. We consider an u.s.c.\ and strongly concave function $\varphi:X\to\R$ with modulus $c>0$ and a Lipschitz continuous function $\psi:X\to\R$ with Lipschitz constant $L_\psi$ with respect to the euclidean norm. Then the maximum of $\varphi$ on $X$ is attained at a unique point $x_0\in X$ i.e.\ $\argmax_X\varphi=\{x_0\}$ and
\begin{equation*}
|\sup_{x\in X}\big(\varphi(x)+\delta\psi(x)\big)-\big(\varphi(x_0)+\delta\psi(x_0)\big)|\leq L_\psi\delta \sqrt{\frac{2\delta M}{c}}\enspace ,
\end{equation*} 
where $M=\sup_{x\in X}\psi(x)-\inf_{x\in X}\psi(x)$.
\end{cor}
\begin{proof}
Define $\Phi(x)=\varphi(x_0)-\varphi(x)$ for $x\in X$ and $\Phi(x)=+\infty$ elsewhere. We have $\Phi(x)\geq 0$ for all $x\in\R^n$ and $\Phi(x_0)=0$. Since $\Phi$ is l.s.c.\ and convex on $\R^n$, then $0\in\partial\Phi(x_0)$. Moreover $\Phi$ is strongly convex with modulus $c$. Then, using Theorem 6.1.2 of~\cite[Chapter VI]{lemarechal2} we have for all $x,x'\in X$
\[
\Phi(x)\geq\Phi(x')+\<s,x-x'>+\frac{c}{2}\|x-x'\|_2^2 \quad \forall s\in\partial\Phi(x')\enspace .
\]
Taking $x'=x_0$ and $s=0$ we obtain for all $x\in X$
\[
\Phi(x)\geq\frac{c}{2}\|x-x_0\|_2^2\enspace ,
\]
which implies
\[
\varphi(x)\leq\varphi(x_0)-\frac{c}{2}\|x-x_0\|_2^2\quad \forall x\in X\enspace .
\]
Using the notations of the previous Lemma, we get easily (see also Proposition 4.32 of~\cite{Bonnans}) for all $x\in F_\varepsilon$, $d(x,F_0)\leq\sqrt{\frac{2\varepsilon}{c}}$, where $\varepsilon=\delta\big(\sup_{x\in X}\psi(x)-\inf_{x\in X}\psi(x)\big)$.
\end{proof}
\begin{rem}\label{erreur-S}
To have an error estimate of the approximation of the matrix $K_{H,h}$ by the matrix $\tilde K_{H,h}$, we apply Lemma~\ref{erreurB} in the case where
\[
\varphi(x)=w_i(x)+z_j(x)\quad\mathrm{and}\quad \psi(x)=H(x,\nabla w_i(x))\enspace ,
\]
for a suitable choice of the finite elements $w_i$ and test functions $z_j$.
Using Assumptions $(H_1)$ and $(H_2)$, we have that, for all $x\in X$, $|\psi(x)|\leq M_f\|\nabla w\|_\infty+M_\ell$, where $\|\nabla w\|_\infty=\|\|\nabla w\|_2\|_\infty$ and $\nabla w=(\nabla w_i)_{1\leq i\leq p}$. We deduce
\[
\sup\psi-\inf\psi\leq2\big(M_f\|\nabla w\|_\infty+M_\ell\big)\enspace .
\]
Moreover $H(\cdot,p)$ and $H(x,\cdot)$ are Lipschitz continuous with Lipschitz constants $L_f\|p\|_2+L_\ell$ and $M_f$ respectively. Hence, $\psi$ is Lipschitz continuous with Lipschitz constant 
\[
L_\psi=L_f\|\nabla w\|_\infty+L_\ell+M_f\|D^2w_i\|_\infty\enspace .
\]
\end{rem}
\subsection{Final estimation of the error of the MFEM}
We now state our main convergence result, which holds for quadratic finite elements and Lipschitz test functions. 
\begin{thm}\label{th-main}
Let $X$ be a compact convex subset of $\R^n$ with non-empty interior and $\hat X=X+\mathrm{B}_2(0,\frac{L}{c})$, where $L>0$, $c>0$. Choose any finite sets of discretization points $\mathcal{T}\subset\R^n$ and $\hat{\mathcal{T}}\subset\R^n$. Let 
\[
\Delta x=\max(\rho_X(\mathcal{T}),\rho_{\hat X}(\hat{\mathcal{T}})).
\]
We make assumptions \mrm{(H1)} and \mrm{(H2)}, and assume that the value function at time $t$, $v^t$, is $c$-semiconvex and Lipschitz continuous with constant $L$ with respect to the euclidean norm, for all $t\geq0$.
Let us choose quadratic finite elements $w_{\hat x}$ of Hessian $c$, centered at the points $\hat x$ of $\hat{\mathcal{T}}$.
Let us choose, as test functions, the Lipschitz finite elements $z_{\hat y}$ with constant $a\geq L$, centered at the points $\hat y$ of $\mathcal{T}$. 
For $t=0,\delta,\ldots, T$, let $v_h^t$ be the approximation of $v^t$ given by the max-plus finite element method implemented with the approximation $K_{H,h}$ of $K_h$ given by~\eqref{appr_K}.
Then, there exists a constant $C_1>0$ such that
\[
\|v_h^T-v^T\|_{\infty}\leq C_1(\delta+\frac{\Delta x}{\delta})\enspace .
\]
When the approximation $K_{H,h}$ is replaced by $\tilde K_{H,h}$, given by~\eqref{e-convenient}, this inequality becomes:
\[
\|v_h^T-v^T\|_{\infty}\leq C_2(\sqrt\delta+\frac{\Delta x}{\delta})\enspace ,
\]
for some constant $C_2>0$.
\end{thm}
\begin{proof}
Let $\calw_h$ and $\calz_h$ denote the complete semimodules of $\rmaxb^{X}$ generated by the families $(w_{\hat x})_{\hat x\in\hat{\mathcal{T}}}$ and $(z_{\hat y})_{\hat y\in\mathcal{T}}$ respectively. We index the elements of $\hat{\mathcal{T}}$ and $\mathcal{T}$ by $\hat x_1,\cdots,\hat x_p$ and $\hat y_1,\cdots,\hat y_q$ respectively. Using Corollary~\ref{error_effect}, we have
\begin{eqnarray*}
\|v_h^T-v^T\|_{\infty}&\leq&(1+\frac T \delta)\Big(\sup_{t\in\bar\tau_\delta}\big(\|P^{-\calz_h}(v^t)-v^t\|_{\infty}+\|P_{\calw_h}(v^t)-v^t\|_\infty\big)\\
&&\quad\quad\quad\quad+\max_{1\leq i\leq p}\|[S^\delta w_{i}]_H-S^\delta w_{i}\|_\infty\Big)\enspace .
\end{eqnarray*}
To estimate the projection error $\|P_{\calw_h}(v^t)-v^t\|_\infty$, we apply Lemma~\ref{proj-p2} for $\hat X_h=\hat{\mathcal{T}}$. We obtain, for $t\in\bar\tau_\delta$, $\|P_{\calw_h}(v^t)-v^t\|_\infty\leq c\operatorname{diam}X\Delta x$. Applying Lemma~\ref{proj-p1} we obtain, for $t\in\bar\tau_\delta$, $\|P^{-\calz_h}(v^t)-v^t\|_{\infty}\leq n(a+L)\Delta x$.
Finally, using Lemma~\ref{appr-StI}, we get
\[
\|v_h^T-v^T\|_{\infty}\leq C_1(\delta+\frac{\Delta x}{\delta})\enspace ,
\]
where 
\[
C_1>(T+1)\max\Big(c\operatorname{diam} X+n(a+L),\frac{1}{2}(L_\ell M_f+L_fM_f(\operatorname{diam} X+\frac{L}{c})+cM_f^2)\Big).
\]
To prove the second inequality, we use Lemma~\ref{complete-error} together with Remark~\ref{erreur-S}. Using the notation of Corollary~\ref{erreurB} and the fact that $\varphi=w_i+z_j$ is c-strongly convex, we have $\sup\psi-\inf\psi\leq 2(M_\ell+M_fc(\operatorname{diam}X+\frac{L}{c}))$ and $L_\psi=L_\ell+cM_f+L_fc(\operatorname{diam}X+\frac{L}{c})$.  We deduce that 
\[
|(\tilde K_{H,h})_{ji}-(K_{H,h})_{ji}|\leq 2\big(L_\ell+cM_f+L_fc(\operatorname{diam}X+\frac{L}{c})\big)\sqrt{\frac{M_\ell}{c}+M_f(\operatorname{diam}X+\frac{L}{c})}\delta\sqrt{\delta}\enspace ,
\]
for $i=1,\cdots,p$ and $j=1,\cdots,q$. Hence, there exists $C2>0$ such that
\[
\|v_h^T-v^T\|_{\infty}\leq C_2(\sqrt\delta+\frac{\Delta x}{\delta})\enspace ,
\]
when $\delta$ is small enough.
\end{proof}
A variant of this theorem, with a stronger assumption, was proved in~\cite{asma}. 
\begin{rem}
When $\mathcal{T}$ is a rectangular grid of step $h>0$, meaning that $\mathcal{T}$ is the intersection of $(\Z h)^n$ with a cartesian product of bounded intervals, we have
\[
\rho_X(\mathcal{T})\leq \sqrt{n}h.
\]
Hence, when $\mathcal{T}$ and $\hat{\mathcal{T}}$ are both rectangular grids of step $h$, we have $\Delta x\leq\sqrt{n}h=O(h)$ in Theorem~\ref{th-main}.
\end{rem}
\section{Numerical results}\label{Numerical results}
This section presents the results of numerical experiments with the MFEM described in Section~\ref{method}. We consider optimal control problems in dimension 1 and 2 whose value functions are known or can be computed by solving the Riccati equation (in the case of linear quadratic problems). 
\subsection{Implementation}
We implemented the MFEM using the max-plus toolbox of Scilab~\cite{toolbox} (in dimension 1) and specific programs written in C (in dimension 2). We used the approximation $\tilde K_{H,h}$ of the matrix $K_h$. The matrix $M_h$ can always be computed analytically. In all the examples below, the Hamiltonian $H$, and so the stiffness matrix $\tilde K_{H,h}$, have been computed analytically. We avoided storing the (full) matrices $M_h$ and $\tilde K_{H,h}$ when the number of discretization points is large.
\subsection{Examples in dimension1}
The next two examples are inspired by those proposed by M. Falcone
in \cite{bardi-capuzzo-dolcetta}.
\begin{exmp}\label{ex-falcone1}
We consider the case where $T=1$, $\phi\equiv 0$, $X=[-1,1]$,
$U=[0,1]$, $\ell(x,u)=x$ and $f(x,u)=-xu$. Assumptions $(H1)$ and $(H2)$ are satisfyied. The optimal choice is to
take $u^*=0$ whenever $x>0$ and to move on the right with maximum
speed ($u^*=1$) whenever $x\leq 0$. For all $t\in [0,T]$, the value function is:
\begin{equation*}
v(x,t)=\begin{cases}
xt & \text{if }  x>0\\
x(1-e^{-t}) & \text{otherwise.}
\end{cases}
\end{equation*}
We choose quadratic finite elements $w_i$ of Hessian $c$ centered at the points of the regular grid $(\Z \Delta x)\cap [-2,2]$ and
Lipschitz finite elements $z_j$ with constant $a\geq 1$ centered at the points of the regular grid $(\Z \Delta x)\cap X$. We
represent in Figure~\ref{fal1} the solution given by our
algorithm in the case where $\delta=0.01$, $\Delta x=0.005$,
$a=1.5$ and $c=1$. We obtain a $L_\infty$-error of order $10^{-2}$.
\begin{figure}[htbp]
\begin{center}
\includegraphics{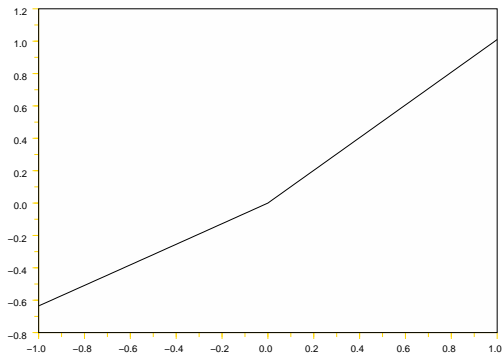}
\caption{Max-plus approximation (Example~\ref{ex-falcone1})}
\label{fal1}
\end{center}
\end{figure}
\end{exmp}
\begin{exmp}\label{ex-falcone2}
We consider the case where $T=1$, $\Phi\equiv 0$, $X=[-1,1]$,
$U=[-1,1]$, $\ell(x,u)=-3(1-|x|)$ and $f(x,u)=u(1-|x|)$. It is clear that $\ell$ and $f$ are bounded and Lipschitz continuous functions. The optimal
choice is to take $u^*=-1$ whenever $x>0$ and $u^*=1$ whenever
$x<0$. Therefore, all the trajectories lie in X. For all $t\in [0,T]$, the value function is:
\begin{equation*}
v(x,t)=-3(1-|x|)(1-e^{-t})
\end{equation*}
We choose quadratic finite elements $w_i$ of Hessian $c$ and
Lipschitz finite elements $z_j$ with constant $a$. We
represent in Figure~\ref{fal2} the solution given by our
algorithm in the case where $\delta=0.02$, $\Delta x=0.01$,
$a=2$ and $c=8$. We obtain a $L_\infty$-error of order $7.66\cdot 10^{-3}$. 
\begin{figure}[htbp]
\begin{center}
\includegraphics{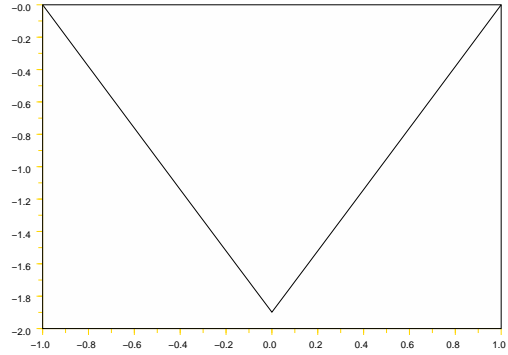}
\caption{Max-plus approximation (Example~\ref{ex-falcone2})}
\label{fal2}
\end{center}
\end{figure}
\end{exmp}

\begin{exmp}[Linear Quadratic Problem]\label{ex-lq}
We consider the case where $U=\R$, $X=\R$,
\[
\ell(x,u)=-\frac{1}{2}(x^2+u^2),
\quad f(x,u)=u,\mrm{ and } \phi\equiv 0
\enspace .
\]
The Hamiltonian is $H(x,p)=-\frac{x^2}{2}+\frac{p^2}{2}$. This problem can be solved analytically. For $x\in X$, the value function at time $t$ is
\[
v(x,t)=-\frac{1}{2}\mathrm{tanh}(t)x^2.
\]
The domain $X$ is unbounded and $\ell$ and $f$ are unbounded and locally Lipschitz continuous. We will restrict $X$ to the set $[-5;5]$ so that $\ell$ and $f$ satisfy Assumptions $(H1)$ and $(H2)$.\\
We choose quadratic finite elements $w_i$ and $z_j$ of Hessian $c=1$, centered at the points of the regular grid $(\Z \Delta x)\cap [-6,6]$.
We represent in Figure~\ref{quadratique} the solution given by
our algorithm in the interval $[-1;1]$ in the case where $T=5$, $\delta=0.5$, $\Delta x=0.05$ and $L=1$. 
We obtain a $L_\infty$-error of $4.54\cdot 10^{-5}$.
\begin{figure}[htbp]
\begin{center}
\includegraphics[scale=0.4]{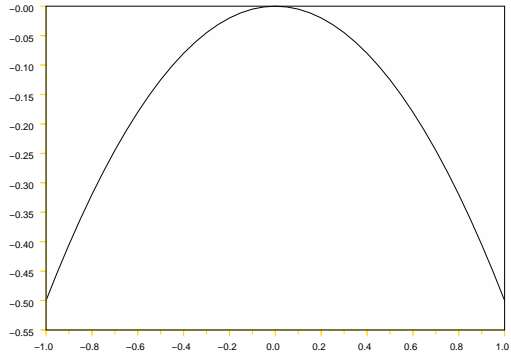}
\caption{Max-plus approximation of a linear quadratic control problem (Example~\ref{ex-lq})}
\label{quadratique}
\end{center}
\end{figure}
\end{exmp}
\begin{exmp}[Distance problem]\label{ex-dist}
We consider the case where $T=1$, $\phi\equiv 0$, $X=[-1,1]$, $U=[-1,1]$,  
\[
\ell(x,u)=\begin{cases}
-1 & \mathrm{if} \quad x \in (-1,1),\\
 0 & \mathrm{if} \quad x \in \{-1,1\} ,
\end{cases}
\quad \mrm{and}\quad
f(x,u)=\begin{cases}
u & \mathrm{if} \quad x \in (-1,1),\\
0 & \mathrm{if} \quad x \in \{-1,1\}.
\end{cases}
\]
Putting $\ell=0$ and $f=0$ on $\partial X$ keeps the trajectories in the domain $X$ but we loose the Lipschitz continuity of $\ell$ and $f$. For $x\in X$, the value function at time $t$ of this problem is
\[
v(x,t)=\max(-t,|x|-1).
\]
Consider first quadratic finite elements $w_i$ and $z_j$ of Hessian $c$, centered at the points of the regular grid $(\Z \Delta x)\cap \big(X+\mathrm{B}_\infty(0,\frac L c)\big)$.
In Figure~\ref{distmauvais}, we represent the solution given by our algorithm in the case where $\delta=0.02$, $\Delta x=0.01$, $c=2$ and $L=1$. 
Since $\Pi^{Z_h^*}$ is a projector on a subsemimodule of the $\rminb$-semimodule of $c$-semiconcave functions, and since the solution is not $c$-semiconcave for any $c$, the error of projection  $\|\Pi^{Z_h^*}(v^t)-v^t\|_\infty$ does not converge to zero when $\Delta x$ goes to zero, which explains the magnitude of the error.
\begin{figure}[htbp]
\begin{center}
\input{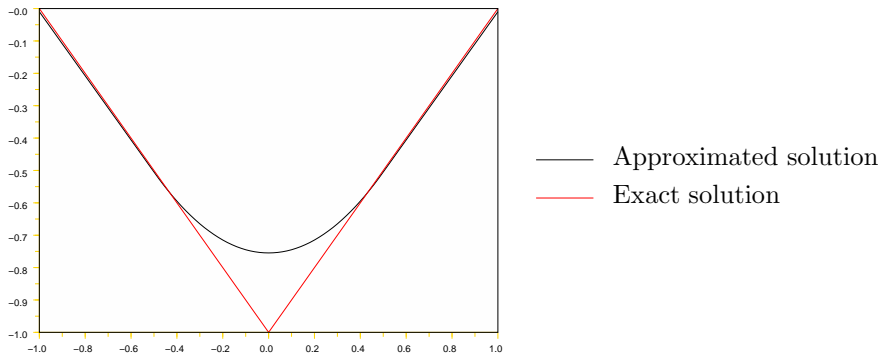}
\caption{A bad choice of test functions for the distance problem (Example~\ref{ex-dist})}
\label{distmauvais}
\end{center}
\end{figure}

To solve this problem, it suffices to replace
the test functions $z_j$ by the Lipschitz finite elements with constant $a\geq 1$, centered at the points of the regular grid
$(\Z \Delta x)\cap [-1,1]$. This is illustrated 
in Figure~\ref{dist} in the case where $\delta=0.02$, $\Delta
x=0.01$, $c=2$ and $a=1.1$. We obtain a $L_\infty$-error of $1.05\cdot 10^{-2}$.
\begin{figure}[htbp]
\begin{center}
\includegraphics{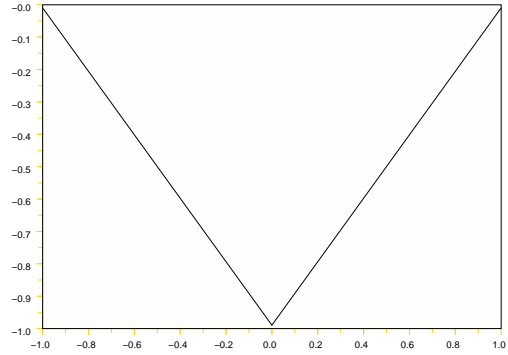}
\caption{A good choice of test functions for the distance problem (Example~\ref{ex-dist})}
\label{dist}
\end{center}
\end{figure}
\end{exmp}
\subsection{Examples in dimension 2}
\begin{exmp}[Linear Quadratic Problem in dimension 2]\label{ex-lq2}
We consider the case where $U=\R^2$, $X=\R^2$, $\phi\equiv 0$,
\[
\ell(x,u)=-\frac{x_1^2+x_2^2}{2}-\frac{u_1^2+u_2^2}{2}\enspace\mrm{ and } \quad f(x,u)=u \enspace .
\]
For $x\in X$, the value functions at time $t$ is
\[
v(x,t)=-\frac{1}{2}\mathrm{tanh}(t)(x_1^2+x_2^2).
\]
As in Example~\ref{ex-lq}, the domain $X$ is unbounded therefore $\ell$ and $f$ do not satisfy Assumptions $(H1)$ and $(H2)$. We will restrict the domain to the set $[-5;5]^2$.\\
We choose quadratic finite elements $w_i$ and $z_j$ of Hessian $c$ centered at the points of the regular grid $\big((\Z \Delta x)\cap [-6,6]\big)^2$.
We represent in Figure~\ref{quadratique2} the solution given by our algorithm in the case where $T=5$, $\delta=0.5$, $\Delta x=0.1$, $c=1$. 
\begin{figure}[htpb]
\begin{center}
\includegraphics{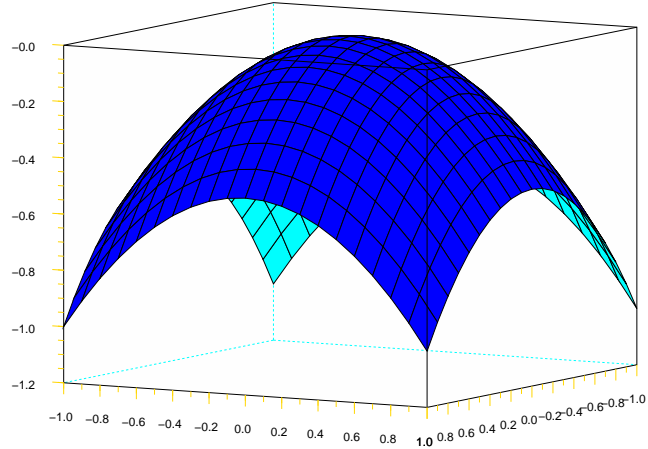}
\caption{Max-plus approximation of a linear quadratic control problem (Example~\ref{ex-lq2})}
\label{quadratique2}
\end{center}
\end{figure}
The $L_{\infty}$-error is $9\cdot 10^{-5}$.
\end{exmp}
\begin{exmp}[Distance problem in dimension 2]\label{ex-dist2}
We consider the case where $T=1$, $\phi\equiv 0$, $X=[-1,1]^2$, $U=[-1,1]^2$,
\[
\ell(x,u)=\begin{cases}
-1 & \mathrm{if} \quad x\in\mathrm{int} X,\\
 0 & \mathrm{if} \quad x\in\partial X,
\end{cases}
\]
\[
f(x,u)=\begin{cases}
u & \mathrm{if} \quad x\in\mathrm{int} X,\\
0 & \mathrm{if} \quad x\in\partial X.
\end{cases}
\]
For $x\in X$, the value function at time $t$ is
\[
v(x,t)=\max\big(-t,\max(|x_1|,|x_2|)-1\big).
\]
We choose quadratic finite elements $w_i$ of Hessian $c$ centered at the points of the regular grid $\big((\Z\Delta x)\cap[-3,3]\big)^2$ and Lipschitz finite elements $z_j$ with constant $a$ centered at the points of the regular grid $\big((\Z\Delta x)\cap[-1,1]\big)^2$. We represent in Figure~\ref{dist2} the solution given by our algorithm in the case where $T=1$, $\delta=0.05$, $\Delta x=0.025$, $a=3$ and $c=1$. The $L_{\infty}$-error is of order $0.05$.
\begin{figure}[thpb]
\begin{center}
\includegraphics{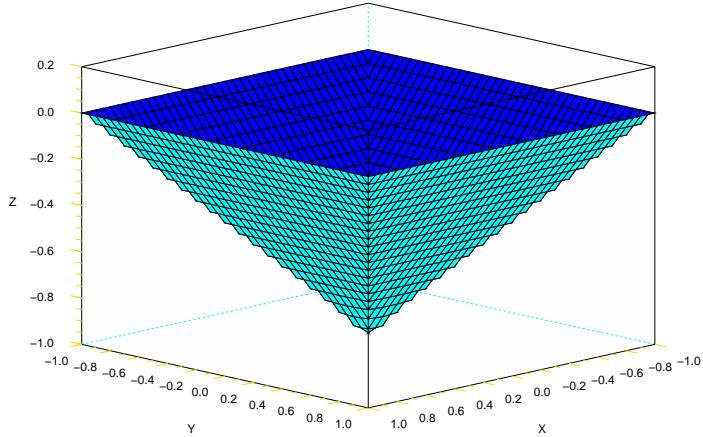}
\caption{Max-plus approximation of the distance problem (Example~\ref{ex-dist2})}
\label{dist2}
\end{center}
\end{figure}
\end{exmp}
\begin{exmp}[Rotating problem]\label{ex-rotation}
We consider here the Mayer problem where $T=1$, $X=\mathrm{B}_2(0,1)$, $U=\{0\}$, $\phi(x)=-\frac 1 2 x_1^2-\frac 3 2 x_2^2$, $\ell(x,u)=0$ and $f(x,u)=(-x_2,x_1)$. For $x\in X$, the value function at time $t$ is
\[
v(x,t)=-\frac 1 2(-x_2\mathrm{sin}(t)+x_1\mathrm{cos}(t))^2-\frac 3 2(x_2\mathrm{cos}(t)+x_1\mathrm{sin}(t))^2.
\]
We choose quadratic finite elements $w_i$ and $z_j$ of Hessians $c_w$ and $c_z$ respectively, centered at the points of the regular grid $\big((\Z \Delta x)\cap [-2,2]\big)^2$.
We represent in Figure~\ref{rotation} the solution given by our algorithm in the case where $\delta=\Delta x=0.05$, $c_w=4$ and $c_z=3$. The  $L_{\infty}$-error is $0.046$.
\begin{figure}[thpb]
\begin{center}
\includegraphics[scale=2]{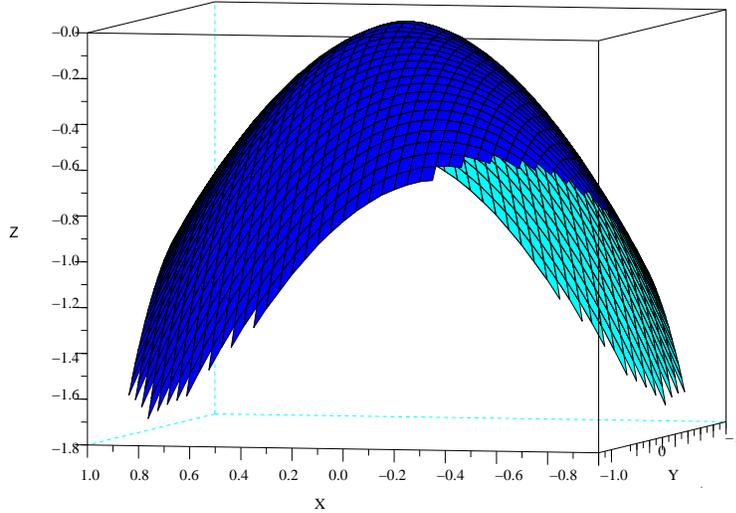}
\caption{Max-plus approximation of the rotating problem (Example~\ref{ex-rotation})}
\label{rotation}
\end{center}
\end{figure}
\end{exmp}
\begin{exmp}\label{ex-riccati}
We consider the case where $U=\R$, $X=\R^2$, $\phi(x)=-x_1^2-2x_2^2$,
\[
\ell(x,u)=-x_1^2-\frac{u^2}{2}\enspace\mrm{ and } \quad f(x,u)=(x_2,u)^T \enspace .
\]
We choose quadratic finite elements $w_i$ and $z_j$ of Hessian $c_w$ and $c_z$ respectively centered at the points of the grids $\big((\Z \Delta x)\cap [-2,2]\big)^2$ and $\big((\Z \Delta x)\cap [-11,11]\big)^2$ respectively.
We represent in Figure~\ref{riccati} the solution given by our algorithm in the case where $T=1$, $\delta=0.05$, $\Delta x=0.025$, $c_w=10$ and $c_z=1$. 
\begin{figure}[htpb]
\begin{center}
\includegraphics{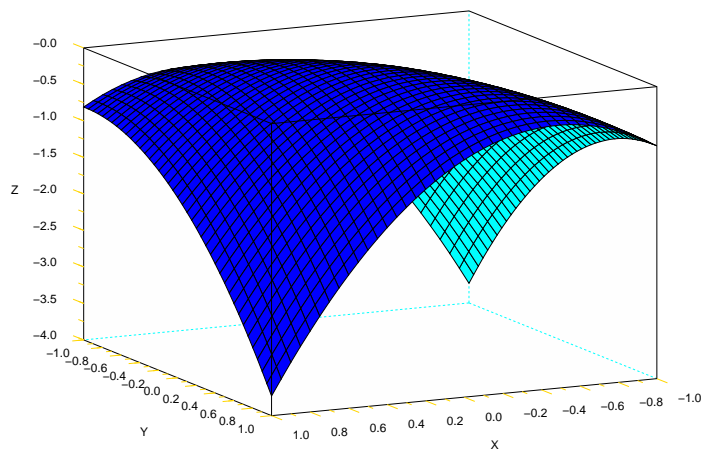}
\caption{Max-plus approximation of the solution of the control problem of Example~\ref{ex-riccati}}
\label{riccati}
\end{center}
\end{figure}
The $L_{\infty}$-error is $0.11$. (We compared the max-plus approximation with the solution of the problem given by the Riccati equation).
\end{exmp}
\subsection{Conclusion}
We have tested our method on examples that fullfill the assumptions of Theorem~\ref{th-main} (see Examples~\ref{ex-falcone1}, \ref{ex-falcone2}, \ref{ex-rotation}) but also on problems that do not fullfill these assumptions. The method is efficient even in the second case. The only difficulty comes from the full character of the matrices $M_h$ and $K_h$, which limits the number of discretization points. To treat higher dimensional examples, we need higher order approximations (when the value function is regular enough). This is the object of a subsequent work.
\bibliographystyle{myalpha} 

\begin{thebibliography}{BCOQ92}
\providecommand{\url}[1]{\texttt{#1}}
\providecommand{\urlprefix}{URL }
  \providecommand{\doi}[1]{Eprint \href{http://dx.doi.org/#1}{doi:#1}}
\providecommand{\arxiv}[2][]{Also \href{http://www.arXiv.org/abs/#2}{arXiv:#2}}

\bibitem[AGL04]{MTNS}
M.~Akian, S.~Gaubert, and A.~Lakhoua.
\newblock A max-plus finite element method for solving finite horizon
  deterministic optimal control problems.
\newblock In \emph{Proceedings of the Sixteenth International Symposium on
  Mathematical Theory of Networks and Systems ({MTNS'04})}. Leuven, Belgium,
  2004.
\newblock And arXiv:math.OC/0404184, April 2004.

\bibitem[Bar94]{barles}
G.~Barles.
\newblock \emph{Solutions de viscosit\'{e} des \'{e}quations de
  {H}amilton-{J}acobi}.
\newblock Springer Verlag, 1994.

\bibitem[BCD97]{bardi-capuzzo-dolcetta}
M.~Bardi and I.~Capuzzo-Dolcetta.
\newblock \emph{Optimal control and viscosity solutions of
  {H}amilton-{J}acobi-{B}ellman equations}.
\newblock Birkha\"user, 1997.

\bibitem[BCOQ92]{baccelli}
F.~Baccelli, G.~Cohen, G.~J. Olsder, and J.-P. Quadrat.
\newblock \emph{Synchronization and linearity~: an algebra for discrete events
  systems}.
\newblock John Wiley \& Sons, New-York, 1992.

\bibitem[BD99]{boue-dupuis}
M.~Bou{\'e} and P.~Dupuis.
\newblock Markov chain approximations for deterministic control problems with
  affine dynamics and quadratic cost in the control.
\newblock \emph{SIAM J. Numer. Anal.}, 36(3):667--695 (electronic), 1999.
\newblock ISSN 0036-1429.

\bibitem[Bir67]{birkhoff}
G.~Birkhoff.
\newblock \emph{Lattice Theory}, volume~25.
\newblock American Mathematical Society, 1967.

\bibitem[BJ72]{blyth}
T.~S. Blyth and M.~F. Janowitz.
\newblock \emph{Residuation theory}.
\newblock Pergamon Press, Oxford, 1972.
\newblock International Series of Monographs in Pure and Applied Mathematics,
  Vol. 102.

\bibitem[BS00]{Bonnans}
J.~F. Bonnans and A.~Shapiro.
\newblock \emph{Perturbation analysis of optimization problems}.
\newblock Springer Series in Operations Research. Springer Verlag, New York,
  2000.

\bibitem[BZ05]{zidani-bokanowski}
O.~Bokanowski and H.~Zidani.
\newblock Anti-dissipative schemes for advection and application to
  hamilton-jacobi-bellman equations.
\newblock \emph{J. Sci. Compt}, to appear 2005.

\bibitem[CD83]{capuzzodolcetta}
I.~Capuzzo~Dolcetta.
\newblock On a discrete approximation of the {H}amilton-{J}acobi equation of
  dynamic programming.
\newblock \emph{Appl. Math. Optim.}, 10(4):367--377, 1983.

\bibitem[CDF89]{capuzzodolcetta-falcone}
I.~Capuzzo-Dolcetta and M.~Falcone.
\newblock Discrete dynamic programming and viscosity solutions of the {B}ellman
  equation.
\newblock \emph{Ann. Inst. H. Poincar\'e Anal. Non Lin\'eaire},
  6(suppl.):161--183, 1989.
\newblock Analyse non lin\'eaire (Perpignan, 1987).

\bibitem[CDI84]{capuzzodolcetta-ishii}
I.~Capuzzo-Dolcetta and H.~Ishii.
\newblock Approximate solutions of the {B}ellman equation of deterministic
  control theory.
\newblock \emph{Appl. Math. Optim.}, 11(2):161--181, 1984.

\bibitem[CFF04]{carlini-falcone-ferretti}
E.~Carlini, M.~Falcone, and R.~Ferretti.
\newblock An efficient algorithm for {H}amilton-{J}acobi equations in high
  dimension.
\newblock \emph{Comput. Vis. Sci.}, 7(1):15--29, 2004.

\bibitem[CG79]{cuning}
R.~Cuninghame-Green.
\newblock \emph{Minimax Algebra}.
\newblock Number 166 in Lecture notes in Economics and Mathematical Systems.
  Springer Verlag, 1979.

\bibitem[CGQ96]{wodes}
G.~Cohen, S.~Gaubert, and J.-P. Quadrat.
\newblock Kernels, images and projections in dioids.
\newblock In \emph{Proceedings of the International Workshop on Discrete Event
  Systems ({WODES'96})}. IEE, Edinburgh, UK, 1996.

\bibitem[CGQ04]{ilade}
G.~Cohen, S.~Gaubert, and J.-P. Quadrat.
\newblock Duality and separation theorem in idempotent semimodules.
\newblock \emph{Linear Algebra and Appl.}, 379:395--422, 2004.
\newblock \doi{10.1016/j.laa.2003.08.010}.
\newblock \arxiv{math.FA/0212294}.

\bibitem[CL84]{crandall-lions}
M.~G. Crandall and P.-L. Lions.
\newblock Two approximations of solutions of {H}amilton-{J}acobi equations.
\newblock \emph{Math. Comp.}, 43(167):1--19, 1984.

\bibitem[CM04]{mceneaney-collins}
G.~Collins and W.~McEneaney.
\newblock Min-plus eigenvector methods for nonlinear {$H\sb \infty$} problems
  with active control.
\newblock In \emph{Optimal control, stabilization and nonsmooth analysis},
  volume 301 of \emph{Lecture Notes in Control and Inform. Sci.}, pages
  101--120. Springer, Berlin, 2004.

\bibitem[CT80]{crandall}
M.~G. Crandall and L.~Tartar.
\newblock Some relations between non expansive and order preserving maps.
\newblock \emph{Proceedings of the AMS}, 78(3):385--390, 1980.

\bibitem[DJLC53]{marie-louise}
M.~Dubreil-Jacotin, L.~Lesieur, and R.~Croisot.
\newblock \emph{Th\'{e}orie des treillis des structures alg\'{e}briques
  ordonn\'{e}es et des treillis g\'{e}om\'{e}triques}.
\newblock Gauthier-Villars, Paris, 1953.

\bibitem[Fal87]{falcone}
M.~Falcone.
\newblock A numerical approach to the infinite horizon problem of deterministic
  control theory.
\newblock \emph{Appl. Math. Optim.}, 15(1):1--13, 1987.
\newblock Corrigenda in {\it Appl. Math. Optim.}, 23:213--214, 1991.

\bibitem[Fat06]{fathi}
A.~Fathi.
\newblock \emph{Weak {KAM} theorem in {L}agrangian dynamics}.
\newblock Cambridge University Press, 2006.
\newblock To appear.

\bibitem[FF94]{falcone-ferretti}
M.~Falcone and R.~Ferretti.
\newblock Discrete time high-order schemes for viscosity solutions of
  {H}amilton-{J}acobi-{B}ellman equations.
\newblock \emph{Numer. Math.}, 67(3):315--344, 1994.
\newblock ISSN 0029-599X.

\bibitem[FG99]{falcone-giorgi}
M.~Falcone and T.~Giorgi.
\newblock An approximation scheme for evolutive {H}amilton-{J}acobi equations.
\newblock In \emph{Stochastic analysis, control, optimization and
  applications}, Systems Control Found. Appl., pages 289--303. Birkh\"auser
  Boston, Boston, MA, 1999.

\bibitem[FM00]{mceneaney}
W.~H. Fleming and W.~M. McEneaney.
\newblock A max-plus-based algorithm for a {H}amilton-{J}acobi-{B}ellman
  equation of nonlinear filtering.
\newblock \emph{SIAM J. Control Optim.}, 38(3):683--710, 2000.
\newblock \doi{10.1137/S0363012998332433}.

\bibitem[FS93]{soner}
W.~H. Fleming and H.~M. Soner.
\newblock \emph{Controlled Markov processes and viscosity solutions}.
\newblock Springer Verlag, New-York, 1993.

\bibitem[GM01]{gondran-minoux}
M.~Gondran and M.~Minoux.
\newblock \emph{Graphes, Dio\"\i des et semi-anneaux}.
\newblock TEC \& DOC, Paris, 2001.

\bibitem[Gon96]{gondran96a}
M.~Gondran.
\newblock Analyse {MINPLUS}.
\newblock \emph{C. R. Acad. Sci. Paris S\'er. I Math.}, 323(4):371--375, 1996.
\newblock ISSN 0764-4442.

\bibitem[GR85]{gonzalez-rofman}
R.~Gonzalez and E.~Rofman.
\newblock On deterministic control problems: an approximation procedure for the
  optimal cost, part {I} and {II}.
\newblock \emph{SIAM J. Control Optim.}, 23(2):242--285, 1985.

\bibitem[HUL93]{lemarechal2}
J.-B. Hiriart-Urruty and C.~Lemar\'echal.
\newblock \emph{Convex analysis and minimization algorithms {I}}.
\newblock Springer Verlag, 1993.

\bibitem[KM88]{kolokltsovmaslov88}
V.~N. Kolokoltsov and V.~P. Maslov.
\newblock The {C}auchy problem for the homogeneous {B}ellman equation.
\newblock \emph{Soviet Math. Dokl.}, 36(2):326--330, 1988.

\bibitem[KM97]{kolokoltsov}
V.~N. Kolokoltsov and V.~P. Maslov.
\newblock \emph{Idempotent analysis and applications}.
\newblock Kluwer Acad. Publisher, 1997.

\bibitem[Lak03]{asma}
A.~Lakhoua.
\newblock \emph{R\'esolution num\'erique de probl\`emes de commande optimale
  d\'eterministe et alg\`ebre max-plus}.
\newblock Rapport de {DEA}, Universit\'e Paris VI, 2003.

\bibitem[Lio82]{lions}
P.-L. Lions.
\newblock \emph{Generalised solutions of {H}amilton-{J}acobi equations}.
\newblock Pitman, 1982.

\bibitem[LMS01]{litvinov}
G.~L. Litvinov, V.~P. Maslov, and G.~B. Shpiz.
\newblock Idempotent functional analysis: an algebraic approach.
\newblock \emph{Math. Notes}, 69(5):696--729, 2001.
\newblock \doi{10.1023/A:1010266012029}.
\newblock \arxiv{math.FA/0009128}.

\bibitem[Mas73]{maslov73}
V.~Maslov.
\newblock \emph{M\'ethodes Operatorielles}.
\newblock Mir, Moscou, 1973.
\newblock French Transl. 1987.

\bibitem[McE02]{mceneaney02}
W.~M. McEneaney.
\newblock Error analysis of a max-plus algorithm for a first-order {HJB}
  equation.
\newblock In \emph{Stochastic theory and control (Lawrence, KS, 2001)}, volume
  280 of \emph{Lecture Notes in Control and Inform. Sci.}, pages 335--351.
  Springer, Berlin, 2002.

\bibitem[McE03]{mceneaney03}
W.~M. McEneaney.
\newblock Max-plus eigenvector representations for solution of nonlinear
  {$H_\infty$} problems: basic concepts.
\newblock \emph{IEEE Trans. Automat. Control}, 48(7):1150--1163, 2003.
\newblock ISSN 0018-9286.

\bibitem[McE04]{mceneaney04}
W.~M. McEneaney.
\newblock Max-plus eigenvector methods for nonlinear {$H_\infty$} problems:
  Error analysis.
\newblock \emph{SIAM J. Control Optim.}, 43(2):379--412 (electronic), 2004.

\bibitem[MH98]{mceneaney-hortona}
W.~M. McEneaney and M.~Horton.
\newblock Max-{P}lus eigenvector representations for nonlinear {$H_\infty$}
  value functions.
\newblock In \emph{Proceedings of the 37th {C}onference on {D}ecision and
  {C}ontrol ({CDC'98})}, pages 3506--3511. IEEE, 1998.

\bibitem[MH99]{mceneaney-hortonb}
W.~M. McEneaney and M.~Horton.
\newblock Computation of max-plus eigenvector representations for nonlinear
  {$H_\infty$} value functions.
\newblock In \emph{Americam {C}ontrol {C}onference}, pages 1400--1404. 1999.

\bibitem[MS92]{maslov92}
V.~P. Maslov and S.~Samborski\u\i, editors.
\newblock \emph{Idempotent analysis}, volume~13 of \emph{Adv. in Sov. Math.}
\newblock AMS, RI, 1992.

\bibitem[OS91]{osher-shu}
S.~Osher and C.-W. Shu.
\newblock High-order essentially nonoscillatory schemes for {H}amilton-{J}acobi
  equations.
\newblock \emph{SIAM J. Numer. Anal.}, 28(4):907--922, 1991.

\bibitem[Plu98]{toolbox}
M.~Plus.
\newblock Documentation of the max-plus toolbox of
  \href{http://www.scilab.org}{Scilab}, 1998.
\newblock Available from ftp://ftp.inria.fr/INRIA/Scilab/contrib/MAXPLUS/.

\bibitem[Roc70]{rockafellar}
R.~T. Rockafellar.
\newblock \emph{Convex analysis}.
\newblock Princeton Mathematical Series, No. 28. Princeton University Press,
  Princeton, N.J., 1970.

\bibitem[SU00]{comp-geom}
J.-R. Sack and J.~Urrutia.
\newblock \emph{Handbook of computational geometry}.
\newblock North-Holland, Amsterdam, 2000.

\end{thebibliography}

\end{document}